\documentclass{article}
\usepackage{amsthm,amsfonts,amssymb,epsfig,graphics,amsmath,fancyheadings}
                                                                                                                                                             
\def\BbbR{\textrm{I${}\!$R}}
\def\R{\textrm{Re\ }}
\def\d{{\delta_*}}
\def\constant{\textrm{constant\ }}

\def\myqed{\vrule height3pt depth2pt width3pt \bigskip}
                                                                                                                                                             
\theoremstyle{plain}
\newtheorem{theorem}{Theorem}
\newtheorem{lemma}{Lemma}
\newtheorem{proposition}{Proposition}
\newtheorem{remark}{Remark}
\newtheorem{definition}{Definition}
                                                                                                                                                             
\title{Stability of undercompressive shock profiles}\bigskip

\author{Peter HOWARD \and Kevin ZUMBRUN}

\begin{document}

\maketitle

\begin{abstract}
Using a simplified pointwise iteration scheme,
we establish nonlinear phase-asymptotic orbital stability of
large-amplitude Lax, undercompressive, overcompressive, 
and mixed under--overcompressive type shock profiles
of strictly parabolic systems of conservation laws
with respect to initial perturbations $|u_0(x)|\le E_0 (1+|x|)^{-3/2}$
in $C^{0+\alpha}$, $E_0$ sufficiently small,
under the necessary conditions of spectral and hyperbolic
stability together with transversality of the connecting profile.
This completes the program initiated by Zumbrun and Howard in \cite{ZH},
extending to the general undercompressive case results obtained for Lax 
and overcompressive shock profiles in \cite{SzX}, \cite{L}, 
\cite{ZH}, \cite{Z.2}, \cite{Ra}, \cite{MZ.1}--\cite{MZ.5}, 
and for special undercompressive profiles in 
\cite{LZ.1}--\cite{LZ.2}, \cite{HZ}. 
In particular, together with spectral results of \cite{Z.6},
our results yield nonlinear stability 
of large-amplitude undercompressive phase-transitional profiles 
near equilibrium of Slemrod's model \cite{Sl.5} for van der Waal gas 
dynamics or elasticity with viscosity--capillarity.
\end{abstract}

\section{Introduction}\label{introduction}

In the series of papers \cite{Z.2}, \cite{MZ.1}--\cite{MZ.5},
Zumbrun and Mascia--Zumbrun,
building on methods introduced in \cite{ZH,HZ}, 
have established nonlinear 
$L^1\cap H^3\to L^p$ (resp. $L^1\cap H^2\to L^p$)
orbital stability, $p>1$,
of large-amplitude Lax-type shock profiles of systems of conservation laws
with viscosity (resp. relaxation), by a simple shock-tracking argument using 
mainly $L^q\to L^p$ bounds on the linearized solution operator.
For precursors of this method, see, e.g.,
\cite{Go.2}, \cite{K.1}--\cite{K.2}, \cite{LZ.1}--\cite{LZ.2}, \cite{ZH},
\cite{HZ}.
See also the alternative arguments carried out in \cite{SzX}, \cite{L}
for small-amplitude Lax-type shock profiles of systems with artificial
(Laplacian) viscosity, 
and in \cite{ZH}, \cite{Ra} for large-amplitude Lax- or overcompressive-type 
profiles of systems with general, possibly degenerate viscosity.

The purpose of the present work is to point out that
a simple pointwise version of the argument of 
\cite{Z.2}, \cite{MZ.1}--\cite{MZ.5} may be applied also to
under-, over-, and mixed under--overcompressive shock profiles 
of strictly parabolic systems, giving a simple and unified
treatment of shock stability independent of the amplitude or
type of the connecting profile, depending only on the necessary, 
Evans-function condition established in \cite{ZH}, \cite{GZ},
\cite{BSZ}, \cite{Z.3}--\cite{Z.4}, equivalent to spectral
and hyperbolic stability plus transversality of the connecting
profile as a solution of the associated traveling-wave ordinary differential
equation (ODE).
In particular, we obtain for the first time
nonlinear stability of general undercompressive
profiles such as arise in phase-transitional gas dynamics and
elasticity \cite{Sl.1}--\cite{Sl.5} or multiphase flow
\cite{AMPZ.1}--\cite{AMPZ.3}, \cite{IMP}, 
extending results obtained for special undercompressive profiles in 
\cite{LZ.1}--\cite{LZ.2}, \cite{HZ}. 

Moreover, the slight additional detail afforded by our pointwise description
is sufficient to give also convergence of the phase shift, or phase-asymptotic
orbital stability (definition recalled below),
which was lacking in \cite{Z.2}, \cite{MZ.1}--\cite{MZ.5}. 
This completes the one-dimensional stability program initiated in \cite{ZH},
at least for systems with strictly parabolic viscosity,
giving a complete characterization of stability analogous
to that obtained by Sattinger in the scalar case \cite{Sat}.
The assumption of strict parabolicity is physically appropriate for 
applications to phase-transitional shock waves, at least in
one dimension; see Remark \ref{strict}.
Indeed, together with the spectral 
analysis of \cite{Z.6},  our results yield one-dimensional 
nonlinear stability of large-amplitude
phase-transitional profiles near equilibrium of
Slemrod's model for van der Waals gas dynamics or elasticity
with viscosity--capillarity,
one of only two large-amplitude stability results that have so far
been obtained for physical models (the other being
stability of profiles of isentropic gamma-law gas dynamics
with $\gamma=1$ \cite{MN}).
Stability of undercompressive profiles for systems with degenerate viscosity
remains an interesting open problem.

\medbreak

Consider a traveling-wave solution
\begin{equation}
 u(x,t)=\bar u(x-st), \qquad 
 \lim_{z\to \pm \infty} \bar u(z)= u_\pm,
\label{profile}
\end{equation}
or ``shock profile'',
of a system of conservation laws 
\begin{equation}
u_t + f(u)_x = (B(u) u_x)_x,
\label{parabolic}
\end{equation}
$x$, $t\in \mathbb{R}$,
$u$, $f\in \mathbb{R}^n$, $B\in \mathbb{R}^{n\times n}$,
corresponding to an ``ideal'', or 
discontinuous shock wave 
\begin{equation}
\label{shock}
 u(x,t)= 
\begin{cases}
u_-& x\le st,\\
u_+& x>st,\\
\end{cases}
\end{equation}
of the associated hyperbolic system
\begin{equation}
 u_t + f(u)_x = 0.
 \label{hyperbolic}
\end{equation}
Without loss of generality (changing to coordinates moving
with the shock), take $s=0$, so that \eqref{profile} becomes
a stationary, or standing-wave solution convenient for stability
analysis.

Following \cite{ZH}, \cite{Z.4}, we make the standard assumptions:

\medskip

({H0}) \quad  $f,B\in C^{2}$.
\medskip

({H1}) \quad $Re \, \sigma(B) > 0$.
\medskip

({H2}) \quad $\sigma (f'(u_\pm))$ real, distinct, and  nonzero.
\medskip

({H3}) \quad $Re\,  \sigma(-ik f'(u_\pm) -k^2 B(u_\pm))< -\theta k^2$ for all real $k$,
some $\theta>0$.
\medskip

({H4}) \quad There exists a solution $\bar u$ of
\eqref{profile}--\eqref{parabolic}, nearby
which the set of all solutions connecting the same values $u_\pm$
forms a smooth manifold $\{\bar u^\delta\}$, $\delta\in \mathcal{U}\subset 
\mathbb{R}^\ell$, $\bar u^0=\bar u$.
\bigskip

\begin{definition}\label{type}
An ideal shock \eqref{shock} is classified as 
{\it undercompressive}, {\it Lax}, 
or {\it overcompressive} type according as $i-n$ 
is less than, equal to, or greater than $1$,
where $i$, denoting the sum of the dimensions $i_-$ and $i_+$ 
of the center--unstable subspace of $df(u_-)$ and the center--stable 
subspace of $df(u_+)$, represents the total number of characteristics
incoming to the shock.  

A viscous profile \eqref{profile} is classified
as {\it pure undercompressive} type if the associated ideal shock
is undercompressive and $\ell=1$, 
{\it pure Lax} type
if the corresponding ideal shock is Lax type and $\ell=i-n$,
and {\it pure overcompressive} type if
if the corresponding ideal shock is overcompressive and $\ell=i-n$,
$\ell$ as in (H4) and $i$ as in Definition \ref{type}.
Otherwise it is classified as {\it mixed under--overcompressive} type;
see \cite{LZ.2}, \cite{ZH}.
\end{definition}

Pure Lax type profiles are the most common type arising in 
standard gas dynamics, while pure over- and undercompressive type 
profiles arise in magnetohydrodynamics (MHD) and phase-transitional models.
Mixed under--overcompressive profiles are also possible,
as described in \cite{LZ.2}, \cite{ZH}, but seldom encountered;
indeed, we do not know a physical example.
In the pure Lax or undercompressive case, 
$\{\bar u^\delta\}=\{\bar u(\cdot-\delta)\}$ is just the
set of all translates of the base profile $\bar u$,
whereas in other cases it involves also deformations of $\bar u$.
For further discussion of existence, structure, and
classification of viscous profiles, see, e.g., 
\cite{LZ.2}, \cite{ZH}, \cite{MZ.2}--\cite{MZ.4}, 
\cite{Z.3}--\cite{Z.5}, and references therein.

\begin{definition}\label{orbital}
The profile $\bar u$ is said to be 
{\it nonlinearly orbitally stable} if
$\tilde u(\cdot,t)$ approaches $\bar u^{\delta(t)}$ as $t\to \infty$,
$\bar u^\delta$ as defined in (H4),
for any solution $\tilde u$ of (\ref{parabolic}) with initial data 
sufficiently close in some norm to the original profile $\bar u$.
If, also, the phase $\delta(t)$ converges to a limiting value $\delta(+\infty)$,
the profile is said to be
{\it nonlinearly phase-asymptotically orbitally stable}. 
\end{definition}
\medbreak

An important result of \cite{ZH}
was the identification of the following {\it stability criterion}
equivalent to $L^1 \to L^p$ 
linearized orbital stability of the profile, $p>1$,
where $D(\lambda)$ as described in \cite{GZ}, \cite{ZH}
denotes the Evans function associated with the linearized
operator $L$ about the profile: 
an analytic function analogous to
the characteristic polynomial of a finite-dimensional operator,
whose zeroes away from the essential spectrum agree in location
and multiplicity with the eigenvalues of $L$. 

\medbreak
($\mathcal{D}$) \quad There exist
precisely $\ell$ zeroes of $D(\cdot)$ in the nonstable half-plane
$\R \lambda \ge 0$, necessarily at the origin $\lambda=0$.
\medbreak

As discussed, e.g., in \cite{ZH}, \cite{Z.3}--\cite{Z.5}, 
under assumptions (H0)--(H4), ($\mathcal{D}$) is equivalent to
(i) {\it strong spectral stability}, $\sigma(L)\subset 
\{\R \lambda \le 0\}\cup \{0\}$, 
(ii) {\it hyperbolic stability} of the associated ideal shock, and 
(iii) {\it transversality} of $\bar u$ as a solution of the connection problem
in the associated traveling-wave ODE,
where hyperbolic stability is defined for Lax and undercompressive
shocks by the Lopatinski condition of \cite{M.1}--\cite{M.3}, \cite{Fre}
and for overcompressive shocks by an analogous long-wave stability
condition \cite{Z.3}.
Here and elsewhere $\sigma$ denotes spectrum of a linearized operator or
matrix.

The stability condition holds always for small-amplitude Lax profiles
\cite{Go.1}, \cite{MN}, \cite{KM}, \cite{KMN}, 
\cite{HuZ}, \cite{PZ}, \cite{FreS},
but may fail for large-amplitude, or nonclassical over- or undercompressive
profiles \cite{AMPZ.1}, \cite{GZ}, \cite{FreZ}, \cite{ZS}, \cite{Z.3}.
It may be readily checked numerically, as described, e.g., in
\cite{Br.1}--\cite{Br.2}, \cite{BrZ}, \cite{BDG}.
It was shown by various techniques 
in \cite{ZH}, \cite{Z.2}, \cite{MZ.1}--\cite{MZ.5}, \cite{Ra}
that the linearized stability condition ($\mathcal{D}$) is also
sufficient for nonlinear orbital stability of Lax or overcompressive profiles
of arbitrary amplitude.  
However, up to now, this result had not been verified in the 
undercompressive case.

In this paper, we present a simple pointwise argument
applicable to shocks of any type, establishing that
($\mathcal{D}$) is sufficient for nonlinear phase-asymptotic
orbital stability.
More precisely, denoting by
\begin{equation}\label{aj}
a_1^\pm<a_2^\pm < \cdots < a_n^\pm
\end{equation} 
the eigenvalues of the limiting convection matrices $A_\pm:= df(u_\pm)$, 
define
\begin{equation}\label{theta}
\theta(x,t):=
\sum_{a_j^-<0}(1+t)^{-1/2}e^{-|x-a_j^-t|^2/Lt}
+ \sum_{a_j^+>0}(1+t)^{-1/2}e^{-|x-a_j^+t|^2/Lt},
\end{equation}
\begin{equation}\label{psi1}
\begin{aligned}
\psi_1(x,t)&:=
\chi(x,t)\sum_{a_j^-<0}
(1+|x|+t)^{-1/2} (1+|x-a_j^-t|)^{-1/2}\\
&\quad+ 
\chi(x,t)\sum_{a_j^+>0}
(1+|x|+t)^{-1/2} (1+|x-a_j^+t|)^{-1/2},\\
\end{aligned}
\end{equation}
and
\begin{equation}\label{psi2}
\begin{aligned}
\psi_2(x,t)&:=
(1-\chi(x,t)) (1+|x-a_1^-t|+t^{1/2})^{-3/2}\\
&\quad +(1-\chi(x,t)) (1+|x-a_n^+t|+t^{1/2})^{-3/2},
\end{aligned}
\end{equation}
where $\chi(x,t)=1$ for $x\in [a_1^-t, a_n^+t]$ and zero otherwise,
and $L>0$ is a sufficiently large constant.

Then, we have the following main theorem.

\begin{theorem}\label{nonlin}
Assuming (H0)--(H4) and the linear stability condition ($\mathcal{D}$),
the profile $\bar u$ is nonlinearly phase-asymptotically
orbitally stable with respect to $C^{0+\alpha}$ initial perturbations 
$|u_0(x)|\le E_0 (1+|x|)^{-3/2}$, $E_0$ sufficiently small.
More precisely, there exist $\delta(\cdot)$ and $\delta(+\infty)$
such that
\begin{equation}
\label{pointwise}
\begin{aligned}
|\tilde u(x,t)-\bar u^{\delta(t)}(x)|&\le C E_0 
(\theta+\psi_1+\psi_2)(x,t),\\
  |\dot \delta (t)|&\le C E_0 (1+t)^{-1},\\
  |\delta(t)-\delta(+\infty)|&\le C E_0(1+t)^{-1/2},
\end{aligned}
\end{equation}
where $\tilde u$ denotes the solution of \eqref{parabolic}
with initial data $\tilde u_0=\bar u+u_0$.
\end{theorem}

In particular, Theorem \ref{nonlin} yields the desired result
of nonlinear stability in the undercompressive or mixed case, effectively
completing the one-dimensional stability analysis initiated in
\cite{ZH}.  

\begin{remark}\label{rate}
Pointwise bound \eqref{pointwise} yields as a corollary
the sharp $L^p$ decay rate 
\begin{equation}
\label{Lp}
|\tilde u(x,t)-\bar u^{\delta(t)}(x)|_{L^p}\le C E_0 
(1+t)^{-\frac{1}{2}(1-\frac{1}{p})}, \quad 1\le p\le \infty.
\end{equation}
\end{remark}

\begin{remark}\label{diffusionwaves}
The profile $\theta$ may be recognized as the superposition of
Gaussian ``approximate diffusion waves'' moving along outgoing characteristic
directions, while the profiles $\psi_1$ and $\psi_2$ respectively
account for nonlinear interactions occurring within the characteristic
cone $[a_1^-t,a_n^+ t]$ and the algebraically decaying tail of the
initial data. 
The profiles $\theta$ and $\psi_j$ 
correspond roughly to the nonlinear diffusion and linearly coupled waves 
$\theta$ and $\eta$\footnote{Defined in \cite{SzX} but not explicitly mentioned 
in either reference, $\eta$ refers to the contribution of quadratic source
terms involving $\theta$ alone, i.e., the second iterate in the nonlinear
iteration through Duhamel's formula.}
 in the more detailed description of the
solution carried out for initial data with the same decay rate
in \cite{L}, \cite{ZH} for the Lax and overcompressive case.
The latter works estimate $\tilde u- \bar u$, which contains
an additional error 
\begin{equation}\label{shocklayer}
\bar u(x)-\bar u^{\delta(t)}(x)\sim 
(\partial \bar u^\delta/\partial \delta)\delta(t) 
=O\big( e^{-\eta |x|}) (1+t)^{-1/2}\big)
\end{equation}
near the shock layer that is not present in our analysis.

A difference of the undercompressive from the Lax or
overcompressive cases is that the time-asymptotic distribution of
mass is no longer determined in a simple way by the mass of the initial data,
making difficult the description of nonlinear diffusion and coupled waves.
We avoid this difficulty by estimating only joint upper bounds and not the 
size or shape of component waves.
Besides undercompressive stability, this yields also considerable
simplification in the
pointwise analysis of Lax and overcompressive profiles.
In particular, we nowhere attempt to identify cancellation
in our estimates of nonlinear interactions, taking into
account only {\it transversality} of interacting signals,
similarly as in (different-family) Glimm interaction estimates for
the hyperbolic case \cite{Gl}.
Compare with the analyses of \cite{SzX, L, Ra} in which
a crucial aspect is to identify cancellation 
in the computation of the linearly coupled wave $\eta$.
Compensating for the lower resolution in our scheme
at the level of diffusion waves is the higher resolution 
afforded by shock tracking, as reflected in the absence of
term \eqref{shocklayer}.
That is, by sufficiently resolving the nondecaying
lowest-order part of the Green function
corresponding to shift in the shock location, we
are able to ignore the details of higher-order parts.
\end{remark}

\begin{remark}\label{multid}
Multidimensional nonlinear $L^1\cap L^\infty \to L^p$ stability, $p\ge 2$, 
has been established using $L^q\to L^p$ resolvent bounds for 
Lax and overcompressive shocks in all dimensions $d\ge 2$ but 
for undercompressive shocks only in dimensions $d\ge 4$; 
see \cite{Z.3}--\cite{Z.5}.
As discussed in \cite{HoZ}, \cite{Z.3}, \cite{HH}, 
stability of undercompressive shock fronts in physical dimensions
$d=2$ and $3$ remains an open question even in the 
(diffusive-dispersive or diffusive-higher order diffusive)
scalar case, for which detailed pointwise Green function bounds
are available.  
An interesting future direction might be to attack this problem by pointwise
methods similar to those of this paper.
\end{remark}

\medskip
\begin{remark}\label{strict}
In one dimension, a change of coordinates reduces Slemrod's
model for van der Waals gas dynamics or elasticity with viscosity--capillarity
to a $2\times 2$ parabolic system with diagonal, strictly parabolic viscosity,
to which our results may be applied; see \cite{Sl.5, Z.6}.
Likewise, in one dimension,
the standard models for imiscible three-phase flow in porous media
may be expressed as 
a $2\times 2$ parabolic system with viscosity that is strictly
parabolic on the interior of the physical state space (the
set of saturations summing to one) and degenerate on the boundary;
see \cite{AMPZ.3}.
Thus, our results generically apply here, too.
However, note that shocks with one state on the boundary
do arise in Riemann problems of physical interest and are
not covered by our theory, nor is the degeneracy of the 
symmetric constant-multiplicity type encountered in real
viscosity models.  This would be an interesting case for
further investigation.
In multi-dimensions, neither of these transformations is possible,
and so a more special analysis of each specific equation
would be required for a stability analysis.
\end{remark}
\medskip

{\bf Plan of the paper.}  In Section \ref{lin}, we recall the linearized
estimates carried out in \cite{ZH}, \cite{MZ.3}.
Assuming certain pointwise convolution estimates, we carry out
in Section \ref{uc} the nonlinear stability analysis of the Lax and
undercompressive case. 
By a slight modification of the argument, we carry out in
Section \ref{oc} the nonlinear stability analysis
of the complementary Lax and overcompressive case.
In Section \ref{mixed}, we establish stability of mixed-type shocks in the
special case $B\equiv \constant$, and in Section \ref{gmixed}
in the general case.
Finally, in Section \ref{estimates}, 
we carry out the deferred convolution estimates,
completing the analysis.
\bigbreak

\section{\bf Linearized estimates}\label{lin}

We begin by recalling the pointwise linearized estimates established
in \cite{ZH}, \cite{MZ.3}, expressed in a streamlined form
convenient for the nonlinear analysis to follow.
Linearizing \eqref{parabolic} 
about a fixed stationary solution $\bar u^{\delta_*}(\cdot)$ gives
\begin{equation}\label{linearized}
u_t=L^{\delta_*}u:= (B^\d u_{x})_{x} -(A^\d  u)_{x},
\end{equation}
where
\begin{equation}\label{coeffs}
A^\d u := df(\bar u^\d(x))u-dB(\bar u^\d(x))(u,\bar u^\d_{x}), \qquad
B^\d := B(\bar u^\d(x)).
\end{equation}

\begin{proposition}[\cite{ZH}, \cite{MZ.3}]\label{greenbounds}
Under the assumptions of Theorem \ref{nonlin}, the Green function
$G(x,t;y)$ associated with the linearized equations \eqref{linearized}
may be decomposed as $G=E+\tilde G$, where 
\begin{equation}\label{E}
E(x,t;y)=\sum_{j=1}^\ell 
\frac{\partial \bar u^\delta(x)}{\partial \delta_j}_{|\delta=\d}e_j(y,t),
\end{equation}
\begin{equation}\label{e}
  e_j(y,t)=\sum_{a_k^{-}>0}
  \left(\textrm{errfn }\left(\frac{y+a_k^{-}t}{\sqrt{4\beta_k^{-}t}}\right)
  -\textrm{errfn }\left(\frac{y-a_k^{-}t}{\sqrt{4\beta_k^{-}t}}\right)\right)
  l_{jk}^{-}(y)
\end{equation}
for $y\le 0$ and symmetrically for $y\ge 0$, with
\begin{equation}\label{ljkbounds}
|l_{jk}^\pm|\le C, \qquad
|(\partial/\partial y)l_{jk}^\pm|\le C\gamma e^{-\eta |y|}, 
\end{equation}
and
\begin{equation}\label{Gbounds}
\begin{aligned}
|\partial_{x,y}^\alpha &\tilde G(x,t;y)|\le \\
& C(t^{-|\alpha|/2}+
|\alpha_y| \gamma e^{-\eta|y|} +|\alpha_x| e^{-\eta|x|})
\Big( \sum_{k=1}^n 
t^{-1/2}e^{-(x-y-a_k^{-} t)^2/Mt} e^{-\eta x^+} \\
&+
\sum_{a_k^{-} > 0, \, a_j^{-} < 0} 
\chi_{\{ |a_k^{-} t|\ge |y| \}}
t^{-1/2} e^{-(x-a_j^{-}(t-|y/a_k^{-}|))^2/Mt}
e^{-\eta x^+}, \\
&+
\sum_{a_k^{-} > 0, \, a_j^{+}> 0} 
\chi_{\{ |a_k^{-} t|\ge |y| \}}
t^{-1/2} e^{-(x-a_j^{+} (t-|y/a_k^{-}|))^2/Mt}
e^{-\eta x^-}\Big), \\
\end{aligned}
\end{equation}
$0\le |\alpha| \le 2$ for $y\le 0$ and symmetrically for $y\ge 0$,
for some $\eta$, $C$, $M>0$, where 
$a_j^\pm$ are as in Theorem \ref{nonlin},  $\beta_k^\pm>0$,
$x^\pm$ denotes the positive/negative
part of $x$,  indicator function $\chi_{\{ |a_k^{-}t|\ge |y| \}}$ is 
$1$ for $|a_k^{-}t|\ge |y|$ and $0$ otherwise,
and $\gamma=1$ in the mixed or undercompressive
case and $0$ in the pure Lax or overcompressive case.
Moreover, all estimates are uniform in the supressed parameter $\d$.
\end{proposition}

\begin{remark}\label{GFremarks}
We will refer to the three differently scaled diffusion kernels
in (\ref{Gbounds}) respectively as the {\it convection} kernel, 
the {\it reflection} kernel, and the {\it transmission} kernel.
The $x$-derivative estimate will only be required in the case 
of mixed-type profiles (see Section 6).  Finally, we recall 
the notation
$$
\textrm{errfn} (z) := \frac{1}{2\pi} \int_{-\infty}^z e^{-\xi^2} d\xi.
$$ 
\end{remark}
\medskip

{\bf Proof of Proposition \ref{greenbounds}.}
This is a restatement of the bounds established in 
\cite{ZH}, \cite{MZ.3} for pure undercompressive, Lax,
or overcompressive type profiles; the same argument
applies also in the mixed under--overcompressive case.
Also, though it was not explicitly stated, uniformity with respect to
$\d$ is a straightforward consequence of the argument.
\myqed

\begin{remark}\label{eboundsrmk}
{}From \eqref{e} and \eqref{ljkbounds}, we obtain by straightforward
calculation (see \cite{MZ.3}) the bounds
\begin{equation}\label{ebounds}
\begin{aligned}
|e_j(y,t)|&\le C\sum_{a_k^->0}
  \left(\textrm{errfn }\left(\frac{y+a_k^{-}t}{\sqrt{4\beta_k^{-}t}}\right)
  -\textrm{errfn }\left(\frac{y-a_k^{-}t}{\sqrt{4\beta_k^{-}t}}\right)\right),\\
|e_j (y,t) &- e_j (y,+\infty)| \le C \textrm{errfn} (\frac{|y|-at}{M\sqrt{t}}), 
\quad \textrm{some}\, a>0 \\
|\partial_t  e_j(y,t)|&\le C t^{-1/2} \sum_{a_k^->0} e^{-|y+a_k^-t|^2/Mt},\\
|\partial_y  e_j(y,t)|&\le C t^{-1/2} \sum_{a_k^->0} e^{-|y+a_k^-t|^2/Mt}\\
&\quad +
C\gamma e^{-\eta|y|}
  \left(\textrm{errfn }\left(\frac{y+a_k^{-}t}{\sqrt{4\beta_k^{-}t}}\right)
  -\textrm{errfn }\left(\frac{y-a_k^{-}t}{\sqrt{4\beta_k^{-}t}}\right)\right),\\
|\partial_y e_j (y,t) &- \partial_y e_j(y,+\infty)|
 \le  C t^{-1/2} \sum_{a_k^->0} e^{-|y+a_k^-t|^2/Mt} \\
|\partial_{yt}  e_j(y,t)|&\le C
(t^{-1}+\gamma t^{-1/2}e^{-\eta|y|}) \sum_{a_k^->0} e^{-|y+a_k^-t|^2/Mt}\\
\end{aligned}
\end{equation}
for $y\le 0$, and symmetrically for $y\ge 0$, where $\gamma$ as above
is one for undercompressive profiles and zero otherwise. 
\end{remark}

\begin{remark}\label{difference}
The main difference between the estimates for the mixed or
undercompressive case $\gamma=1$
and the pure Lax or overcompressive case $\gamma=0$
is the presence of slower-decaying $e^{-\theta|y|}$ terms in derivative
estimates for $e$, $\tilde G$.  
As discussed in \cite{LZ.1}--\cite{LZ.2}, \cite{ZH}, \cite{Z.6}, \cite{Z.2},
these are not only technical artifacts, but reflect real differences
in behavior in the undercompressive case:
specifically, that shock dynamics are not governed solely by conservation 
of mass, as in the Lax or overcompressive case,
but by more complicated dynamics of front interaction as indicated by
rapidly decaying modes $\sim e^{-\theta|y|}$.
\end{remark}

\section{Stability of Lax or undercompressive profiles}\label{uc}

We now carry out the proof of Theorem \ref{nonlin} in the Lax
or undercompressive case, which may be treated by 
a particularly simple
argument. In these cases $\ell=1$, and $\bar u^\delta=\bar u(x-\delta)$,
so that we may conveniently work with the ``centered'' perturbation
variable
\begin{equation}\label{pert}
u(x,t):=\tilde u(x+\delta(t), t)-\bar u(x),
\end{equation}
for which \eqref{parabolic} becomes
\begin{equation}
\label{perteq}
u_t-Lu=Q(u,u_x)_x+\dot \delta (t)(\bar u_x + u_x),
\end{equation}
$L:=L^0$, where
\begin{equation}\label{Q}
\begin{aligned}
Q(u,u_x)&=\mathcal{O}(|u|^2 + |u||u_x|)\\
Q(u,u_x)_x&=\mathcal{O}(|u||u_x|+|u_x|^2 + |u||u_{xx}|)\\
\end{aligned}
\end{equation}
so long as $|u|$ remains bounded.

Recalling the standard fact that $\bar u'$ is a stationary
solution of the linearized equations \eqref{linearized}, 
$L\bar u'=0$, or
$$
\int^\infty_{-\infty}G(x,t;y)\bar u_x(y)dy=e^{Lt}\bar u_x(x)
=\bar u'(x),
$$
we have by Duhamel's principle:
$$
\begin{array}{l}
  \displaystyle{
  u(x,t)=\int^\infty_{-\infty}G(x,t;y)u_0(y)\,dy } \\
  \displaystyle{\qquad
  +\int^t_0 \int^\infty_{-\infty} G_y(x,t-s;y)
  (Q(u,u_x)+\dot \delta u ) (y,s)\,dy\,ds + \delta (t)\bar u'(x).}
\end{array}
$$
Defining 
\begin{equation}
 \begin{array}{l}
  \displaystyle{
  \delta (t)=-\int^\infty_{-\infty}e(y,t) u_0(y)\,dy }\\
  \displaystyle{\qquad
  +\int^t_0\int^{+\infty}_{-\infty} e_{y}(y,t-s)(Q(u,u_x)+
  \dot \delta\, u)(y,s) dy ds, }
  \end{array}
 \label{delta}
\end{equation}
following \cite{ZH}, \cite{Z.2}, \cite{MZ.1}--\cite{MZ.2},
where $e$ is defined as in \eqref{e} (that is, $e=\sum_j e_j$),
and recalling the decomposition $G=E+\tilde G$,
we obtain finally the {\it reduced equations}:
\begin{equation}
\begin{array}{l}
 \displaystyle{
  u(x,t)=\int^\infty_{-\infty} \tilde G(x,t;y)u_0(y)\,dy }\\
 \displaystyle{\qquad
  -\int^t_0\int^\infty_{-\infty}\tilde G_y(x,t-s;y)(Q(u,u_x)+
  \dot \delta u)(y,s) dy \, ds, }
\end{array}
\label{u}
\end{equation}
and, differentiating (\ref{delta}) with respect to $t$,
and observing that 
$e_y (y,s)\rightharpoondown 0$ as $s \to 0$, as the difference of 
approaching heat kernels:
\begin{equation}
 \begin{array}{l}
 \displaystyle{
  \dot \delta (t)=-\int^\infty_{-\infty}e_t(y,t) u_0(y)\,dy }\\
 \displaystyle{\qquad
  +\int^t_0\int^{+\infty}_{-\infty} e_{yt}(y,t-s)(Q(u,u_x)+
  \dot \delta u)(y,s)\,dy\,ds. }
 \end{array}
\label{deltadot}
\end{equation}
\medskip

We shall make use of the following three technical lemmas,
the proofs of which are given in Section \ref{estimates}.

\begin{lemma} 
[Short-time theory]
\label{shorttime}
Under the assumptions of Theorem \ref{nonlin},
for data $u_0\in C^{0+\alpha}(x)$, equations \eqref{delta}--\eqref{u} 
(alternatively, \eqref{ocdeltaalt}, \eqref{ocu}
of the following section) admit a unique local solution
$u\in C^{0+\alpha}(x)\cap C^{0+\alpha/2}(t)$, $\delta \in C^{1+\alpha/2}(t)$, 
extending so long as $|u|_{C^{0+\alpha}}$ 
remains bounded.  Moreover, on this domain, 
$\sup_z |u|(\theta+\psi_1+\psi_2)^{-1}(z,\cdot)$ remains continuous
so long as it and $|\dot \delta(1+t)|$ are uniformly bounded 
and, for $t\ge \tau>0$ sufficiently small, 
\begin{equation}\label{deriv}
\sup_{z}|u_x|(\theta+\psi_1+\psi_2)^{-1}(z,t)
\le C\tau^{-1/2}\sup_{z}|u|
(\theta+\psi_1+\psi_2)^{-1}(z,t-\tau).
\end{equation}
\end{lemma}

\begin{lemma}[Linear 
estimates]\label{iniconvolutions}
Under the assumptions of Theorem \ref{nonlin},
\begin{equation}\label{iniconeq}
\begin{aligned}
\int_{-\infty}^{+\infty}|\tilde G(x,t;y)|(1+|y|)^{-3/2}\, dy
&\le C(\theta+\psi_1+\psi_2)(x,t),\\
\int_{-\infty}^{+\infty}|e_t(y,t)|(1+|y|)^{-3/2}\, dy
&\le C(1+t)^{-3/2},\\
\int_{-\infty}^{+\infty}|e(y,t)|(1+|y|)^{-3/2}\, dy
&\le C,\\
 \int^{+\infty}_{-\infty} |e(y,t)-e(y,+\infty)| (1+|y|)^{-3/2}\, dy 
&\le C(1+t)^{-1/2},\\
\end{aligned}
\end{equation}
for $0\le t\le +\infty$, any $a$, $M>0$, for some $C>0$, where $\tilde G$ and 
$e$ are defined as in Proposition \ref{greenbounds}.
\end{lemma}

\begin{lemma}[Nonlinear 
estimates]\label{convolutions}
Under the assumptions of Theorem \ref{nonlin},
\begin{equation}\label{coneq}
\begin{aligned}
\int_0^t\int_{-\infty}^{+\infty}|\tilde G_y(x,t-s;y)|\Psi(y,s)\, dy ds
&\le C(\theta+\psi_1+\psi_2)(x,t),\\
\int_0^t\int_{-\infty}^{+\infty}|e_{yt}(y,t-s)|\Psi(y,s)\, dy ds
&\le C(1+t)^{-1},\\
\int^{+\infty}_0\int^{+\infty}_{-\infty} 
|e_y(y,+\infty)|\Psi(y,s) \,dy\,ds &\le C\gamma, \\
\int_0^t\int_{-\infty}^{+\infty}
|e_y(y,t-s)- e_y(y,+\infty)| \Psi(y,s)\, dyds
&\le C(1+t)^{-1/2},\\
\int_t^{+\infty} \int_{-\infty}^{+\infty}|e_y(y,t-s)|\Psi(y,s)\, dy
&\le C(1+t)^{-1/2},\\
\end{aligned}
\end{equation}
for
\begin{equation}\label{source}
\begin{aligned}
\Psi(y,s)&:=
(1+s)^{1/2}s^{-1/2}(\theta + \psi_1+\psi_2)^2(y,s)\\
&\qquad +
(1+s)^{-1} (\theta+\psi_1+\psi_2)(y,s)
\end{aligned}
\end{equation}
\end{lemma}

\medskip

{\bf Proof of Theorem \ref{nonlin}, Lax or undercompressive case.} 
Define
\begin{equation}
\label{zeta2}
 \zeta(t):= \sup_{y, 0\le s \le t}
 \Big( |u|(\theta+\psi_1+\psi_2)^{-1}(y,t)
 + |\dot \delta (s)|(1+s) \Big).
\end{equation}
We shall establish:

{\it Claim.} For all $t\ge 0$ for which a solution exists with
$\zeta$ uniformly bounded by some fixed, sufficiently small constant,
there holds
\begin{equation}
\label{claim}
\zeta(t) \leq C_2(E_0 + \zeta(t)^2).
\end{equation}
\medskip

{}From this result, provided $E_0 < 1/4C_2$, 
we have that $\zeta(t)\le 2C_2E_0$ implies
$\zeta(t)< 2C_2E_0$, and so we may conclude 
by continuous induction that
 \begin{equation}
 \label{bd}
  \zeta(t) < 2C_2E_0
 \end{equation}
for all $t\geq 0$.
(By Lemma \ref{shorttime}, $u\in C^1$ exists and $\zeta$ remains
continuous so long as $\zeta$ remains bounded by some uniform constant,
hence \eqref{bd} is an open condition.)
Thus, it remains only to establish the claim above.
\medskip

{\it Proof of Claim.}
We must show that $u(\theta+\psi_1+\psi_2)^{-1}$ and
$|\dot \delta(s)|(1+s)$
are each bounded by
$C(E_0 + \zeta(t)^2)$,
for some $C>0$, all $0\le s\le t$, so long as $\zeta$ remains
sufficiently small.

By \eqref{zeta2}, combined with \eqref{deriv}, we have
for $t\ge 1$ that
\begin{equation}\label{layerA}
\begin{aligned}
|u_x(x,t)|&\le C\zeta(t-1)(\theta +\psi_1+\psi_2)(x,t-1)\\
&\le C_2\zeta (t) (\theta +\psi_1+\psi_2)(x,t)\\
\end{aligned}
\end{equation}
and for $0\le t\le 1$ that
\begin{equation}\label{layer2}
\begin{aligned}
|u_x(x,t)|&\le Ct^{-1/2}\zeta(0)(\theta +\psi_1+\psi_2)(x,0)\\
&\le C_2\zeta (t)t^{-1/2} (\theta +\psi_1+\psi_2)(x,t).\\
\end{aligned}
\end{equation}
Combining these estimates, and recalling definition \eqref{zeta2},
we obtain for all $t\ge 0$ and some $C>0$ that
\begin{equation}\label{ubounds}
\begin{aligned}
|\dot \delta(t)|&\le \zeta(t)(1+t)^{-1},\\
|u(x,t)| &\le \zeta (t)(\theta +\psi_1+\psi_2)(x,t),\\
|u_x(x,t)| &\le C\zeta (t)(1+t)^{1/2} t^{-1/2} (\theta +\psi_1+\psi_2)(x,t)\\
\end{aligned}
\end{equation}
and therefore
\begin{equation}\label{Nbounds}
\begin{aligned}
|(Q(u,u_x)+ \dot \delta u)(y,s)|&\le
C\Psi(y,s)
\end{aligned}
\end{equation}
with $\Psi$ as defined in \eqref{source}.

Combining \eqref{Nbounds} with representations
(\ref{u})--(\ref{deltadot}) and
applying Lemmas \ref{iniconvolutions} and \ref{convolutions}, we obtain
$$
 \begin{aligned}
  |u(x,t)| &\le
  \int^\infty_{-\infty} |\tilde G(x,t;y)| |u_0(y)|\,dy
   \\
 &\qquad +\int^t_0
  \int^\infty_{-\infty}|\tilde G_y(x,t-s;y)||(Q(u,u_x)+
  \dot \delta u)(y,s)| dy \, ds \\
  & \le
  E_0 \int^\infty_{-\infty} |\tilde G(x,t;y)|(1+|y|)^{-3/2}\,dy
   \\
 &\quad +
C\zeta(t)^2 \int^t_0
  \int^\infty_{-\infty}|\tilde G_y(x,t-s;y)|
\Psi(y,s) dy \, ds \\
&\le
C(E_0+\zeta(t)^2)(\theta + \psi_1+\psi_2)(x,t)
\end{aligned}
$$
and, similarly,
$$
\begin{aligned}
 |\dot \delta(t)| &\le \int^\infty_{-\infty}|e_t(y,t)|
  |u_0(y)|\,dy \\
  &\qquad +\int^t_0\int^{+\infty}_{-\infty} |e_{yt}(y,t-s)||(Q(u,u_x)+
 \dot \delta u)(y,s)|\,dy\,ds\\
&\le \int^\infty_{-\infty}|e_t(y,t)|(1+|y|)^{-3/2}\, dy
+ \int^t_0\int^{+\infty}_{-\infty} |e_{yt}(y,t-s)|\Psi(y,s) \,dy\,ds\\
&\le C(E_0+\zeta(t)^2)(1+t)^{-1}.\\
\end{aligned}
$$
Dividing by $(\theta+\psi_1+\psi_2)(x,t)$ and $(1+t)^{-1}$,
respectively, we obtain \eqref{claim} as claimed.
\medbreak

{}From \eqref{claim}, we obtain global existence, with
$\zeta(t)\le 2CE_0$.
{}From the latter bound and the definition of $\zeta$
in \eqref{zeta2} we obtain the first two bounds of \eqref{pointwise}.
It remains to establish the third bound, expressing convergence
of phase $\delta$ to a limiting value $\delta(+\infty)$.

By Lemmas \ref{iniconvolutions}--\ref{convolutions}
together with the previously
obtained bounds \eqref{Nbounds} and $\zeta\le CE_0$, 
and the definition \eqref{zeta2} of $\zeta$, the formal limit
$$
 \begin{aligned}
 \delta(+\infty)&:=  \int^\infty_{-\infty}
e(y,+\infty) u_0(y)\,dy  \\
 &\qquad
 +\int^{+\infty}_0\int^{+\infty}_{-\infty} 
e_y(y,+\infty)(Q(u,u_x)+
  \dot \delta u)(y,s)\,dy\,ds \\
&\le 
 \int^\infty_{-\infty}
E_0 |e(y,+\infty)|(1+|y|)^{-3/2} \,dy  \\
 &\qquad
 +\int^{+\infty}_0\int^{+\infty}_{-\infty} 
CE_0 |e_y(y,+\infty)| \Psi(y,s) \,dy\,ds \\
&\le CE_0
\end{aligned}
$$
is well-defined, as the sum of absolutely convergent integrals.

Applying Lemmas \ref{iniconvolutions}--\ref{convolutions} a final time,
we obtain
$$
 \begin{aligned}
 |\delta(t)-\delta(+\infty)| &\le \int^\infty_{-\infty}
|e(y,t)-e(y,+\infty)| | u_0(y)| \,dy  \\
 &\qquad
 +\int^t_0\int^{+\infty}_{-\infty} |e_{y}(y,t-s)-e_y(y,+\infty)|
|(Q(u,u_x)+
  \dot \delta u)(y,s)|\,dy\,ds \\
 &\qquad
 +\int_t^{+\infty}\int^{+\infty}_{-\infty} |e_y(y,+\infty)|
|(Q(u,u_x)+
  \dot \delta u)(y,s)|\,dy\,ds \\
&\le \int^\infty_{-\infty}
|e(y,t)-e(y,+\infty)| (1+|y|)^{-3/2} \,dy  \\
 &\qquad
 +\int^t_0\int^{+\infty}_{-\infty} |e_{y}(y,t-s)-e_y(y,+\infty)|
\Psi(y,s) \,dy\,ds \\
 &\qquad
 +\int_t^{+\infty}\int^{+\infty}_{-\infty} |e_y(y,+\infty)|
\Psi(y,s)\,dy\,ds \\
&\le CE_0(1+t)^{-1/2},
\end{aligned}
$$
establishing the remaining bound and completing the proof.
\myqed

\bigbreak
\section{Overcompressive profiles.}\label{oc}

We may treat the overcompressive case by a slight modification
of the argument of Section \ref{uc}, which applies also to
the Lax case.
As the Lax and overcompressive case have already been treated by different
means in \cite{ZH}, \cite{Ra}
we shall only sketch the changes necessary for 
the argument, omitting most details.

\medbreak
{\bf Modified equations.}
In the overcompressive case, $\ell>1$, $\bar u^\delta$
consists not only of translates of $\bar u$, but also of
orbits distinct from $\bar u$.
In particular, the different representatives are not all
derived from a group action, and so we cannot use a centering
transformation as in \eqref{pert}, which consists of the group
operations
$$
T_{\delta}(\tilde u-T_{-\delta}\bar u)=T_{\delta}\tilde u -\bar u,
$$
where $T_\alpha v(x,t):=v(x+\alpha,t)$ denotes translation
in $x$.
Accordingly, we work with the primitive variable
\begin{equation}\label{ocpert}
u(x,t):=\tilde u(x, t)-\bar u^{\delta(t)}(x)
\end{equation}
and center the {\it equations} instead, about some strategically
chosen $\d$, obtaining in place of \eqref{perteq}
the modified perturbation equation
\begin{equation}
\label{ocperteq}
u_t-L^\d u=Q^\d (u,u_x)_x+\dot \delta (t)
(\partial \bar u^\delta/\partial \delta)_{|\d}
+ R^\d (\delta, u, u_x )_x +
S^\d (\delta, \delta_t),
\end{equation}
where $Q^\d $ is as in \eqref{Q} and
\begin{equation}\label{R}
\begin{aligned}
R^\d &=
\Big(A(\bar u^\d(x))-A(\bar u^{\delta(t)}(x))\Big)u
+
\Big(B(\bar u^\d(x))-B(\bar u^{\delta(t)}(x))\Big)u_x\\
&=\mathcal{O}( e^{-\eta|y|}|\delta-\d|(|u|+|u_x|) )
\end{aligned}
\end{equation}
and
\begin{equation}\label{S}
\begin{aligned}
S^\d &= 
\dot \delta \Big( (\partial \bar u^\delta/\partial \delta)_{|\delta(t)}-
(\partial \bar u^\delta/\partial \delta)_{|\d}
\Big)
=\mathcal{O}( e^{-\eta|y|}|\dot \delta|| \delta-\d| )
\end{aligned}
\end{equation}
account for ``centering errors''; see, e.g., \cite{GoM},
\cite{HoZ}, \cite{Z.4} for related computations.

Defining 
\begin{equation}
\label{ocdelta}
 \begin{aligned}
  \delta (t)&:= \d -\int^\infty_{-\infty}e(y,t) u_0^\d(y)\,dy 
  -\int^t_0\int^{+\infty}_{-\infty} e(y,t-s)
S^\d(\delta, \delta_t)(y,s) dy ds 
\\
  &\qquad
  +\int^t_0\int^{+\infty}_{-\infty} e_{y}(y,t-s)
(Q^\d(u,u_x)+ R^\d(\delta, u,u_x ))
(y,s) dy ds, \\
u_0^\d&:= \tilde u_0 -\bar u^\d\\
\end{aligned}
\end{equation}
by analogy with \eqref{delta}, we obtain reduced equations
\begin{equation}
\label{ocu}
\begin{aligned}
  u(x,t)&=\int^\infty_{-\infty} \tilde G(x,t;y)u_0^\d(y)\,dy \\
&\qquad -
  \int^t_0\int^\infty_{-\infty}\tilde G(x,t-s;y)
  S^\d(\delta, \delta_t))(y,s) dy \, ds \\
 &\qquad
  -\int^t_0\int^\infty_{-\infty}\tilde G_y(x,t-s;y)(Q^\d(u,u_x)+
  R^\d(\delta, u,u_x))(y,s) dy \, ds
\end{aligned}
\end{equation}
and
\begin{equation}
\label{ocdeltadot}
 \begin{aligned}
  \dot \delta (t)&=-\int^\infty_{-\infty}e_t(y,t) u_0^\d(y)\,dy \\
&\qquad
  - \int^t_0\int^{+\infty}_{-\infty} e_{t}(y,t-s)S^\d(\delta,\delta_t)
  (y,s)\,dy\,ds\\
 &\qquad
  +\int^t_0\int^{+\infty}_{-\infty} e_{yt}(y,t-s)(Q^\d(u,u_x)+
  R^\d(\delta, u,u_x))(y,s)\,dy\,ds. 
 \end{aligned}
\end{equation}
\medskip

\medbreak
{\bf Asymptotic shock location.} 
An important feature of the Lax or overcompressive case is that
the stability criterion ($\mathcal{D}$) in these cases
implies that the (nonlinear) $L^1$-asymptotic state of the perturbed
shock must be formally determined by conservation of mass, in the sense that
relation
\begin{equation}\label{massdist}
\begin{aligned}
\int^{+\infty}_{-\infty} (\tilde u_0 (y)-\bar u(y)) dy &= 
\int (\bar u^{\delta_{+\infty}} (y) - \bar u (y) ) dy \\
&\qquad  + 
\sum_{a_j^+ >0}  m^+_j r^+_j 
+\sum_{a_j^- < 0}  m^-_j r^-_j
\end{aligned}
\end{equation}
is full rank, hence uniquely soluble for $m^\pm_j$, $\delta_{+\infty}$
for $\int_{-\infty}^{+\infty}(\tilde u-\bar u)(y)dy$ sufficiently small,
where $a_j^\pm$ and $r^\pm_j$ denote eigenvalues and right eigenvectors
of $A_\pm=f'(u_\pm)$,
$m^\pm_j$ denotes asymptotic mass in the $j$th outgoing 
characteristic field at $\pm \infty$, and $\delta_{+\infty}$ denotes
the asymptotic shock location, or, equivalently,
\begin{equation}\label{zeromass1}
\int^{+\infty}_{-\infty} 
(\tilde u_0 -\bar u^{\delta_{+\infty}})(y)\, dy=
\sum_{a_j^\pm \gtrless 0}  m^\pm_j r^\pm_j
\end{equation}
for $\delta_{+\infty}=\mathcal{O}(E_0)$, the latter estimate
a consequence of full rank;
for further discussion, see \cite{LZ.1}--\cite{LZ.2}, 
\cite{ZH}, \cite{Z.3}, and references therein.
Centering about $\d=\delta_{+\infty}$,
we may thus arrange that
\begin{equation}\label{zeromass}
\int^{+\infty}_{-\infty} u_0^\d (y)=
\int^{+\infty}_{-\infty} 
(\tilde u_0 -\bar u^{\d})(y)\, dy=
\sum_{a_j^\pm \gtrless 0}  m^\pm_j r^\pm_j,
\end{equation}
while maintaining our assumptions on initial perturbation $u_0^\d$.
In these coordinates, we may expect that
$|\delta(t)-\d|$ decays to zero.

A second consequence of \eqref{massdist}, this time at the 
linearized level, is that $e(y,+\infty)$ 
(constant in the Lax or overcompressive case; see Proposition
\ref{greenbounds}) must be orthogonal to all ``outgoing modes''
$r_j^+$, $a_j^+>0$ and and $r_j^-$, $a_j^-<0$, hence
with choice of coordinates \eqref{zeromass} we have
\begin{equation}\label{zed0}
\int_{-\infty}^{+\infty} e(y,+\infty)u_0^\d(y) \, dy=  0.
\end{equation}

A final consequence of \eqref{massdist} is that 
the map 
\begin{equation}\label{param}
\delta \to  
\hat \delta:=\int_{-\infty}^{+\infty}\Pi (\bar u^\delta- \bar u)(y)\, dy
\in \BbbR^\ell
\end{equation}
must be invertible, where $\Pi\in \ell \times n$ is any (constant) 
full rank matrix with rows orthogonal to outgoing modes $r_j^\pm$.
Reparametrizing by $\delta=\hat \delta$, we may thus arrange
that $\delta =\int_{-\infty}^{+\infty}\Pi (\bar u^\delta- \bar u)(y)\, dy$,
and thus
\begin{equation}\label{Id}
\int_{-\infty}^{+\infty}\Pi (\partial \bar u^\delta/\delta)(y)\, dy
\equiv I_\ell
\end{equation}
for $\delta$ in a neighborhood of the origin, so that
$\int_{-\infty}^{+\infty} \Pi S^\d(y,s)\, dy \equiv 0$ for all $s$, and
therefore
\begin{equation}\label{zed}
\int_{-\infty}^{+\infty}e(y,+\infty)S^\d(y,s)\, dy \equiv 0.
\end{equation}

Combining \eqref{zed0} and \eqref{zed}, we obtain
the alternative representation 
\begin{equation}
\label{ocdeltaalt}
 \begin{aligned}
  \delta (t)&= \d -\int^\infty_{-\infty}
(e(y,t)- e(y,+\infty)) u_0(y)\,dy \\
&\qquad
  -\int^{t}_0\int^{+\infty}_{-\infty} 
(e(y,t-s)- e(y,+\infty))
S^\d(\delta, \delta_t)(y,s) dy ds 
\\
  &\qquad
  +\int^{t}_0\int^{+\infty}_{-\infty} e_{y}(y,t-s)
(Q^\d(u,u_x)+ R^\d(\delta, u,u_x))
(y,s) dy ds,\\
  \end{aligned}
\end{equation}
from which we may observe decay in $|\delta-\d|$
{\it without a priori knowledge of the global behavior of $u$}.

\medbreak
{\bf Stability argument.}
Working within the framework of equations \eqref{ocu}, 
\eqref{ocdeltadot}, and \eqref{ocdeltaalt},
we can carry out the proof of Theorem \ref{nonlin}
for the Lax or overcompressive case by essentially the same argument
presented for the Lax and undercompressive case in the previous section
using Lemmas \ref{shorttime}--\ref{convolutions} together
with the following Lemma proved in Section \ref{estimates}.

\begin{lemma}[Auxiliary estimates]\label{occonvolutions}
Under the assumptions of Theorem \ref{nonlin},
\begin{equation}\label{occoneq1}
\begin{aligned}
\int_0^t\int_{-\infty}^{+\infty} |\tilde G_y(x,t-s;y)|
\Phi_1(y,s) \, dy ds
&\le C(\theta+\psi_1+\psi_2)(x,t),\\
\int_0^t\int_{-\infty}^{+\infty} |e_{yt}(y,t-s)|
\Phi_1(y,s) \, dy ds
&\le C(1+t)^{-1},\\
\int_0^t\int_{-\infty}^{+\infty} |e_{y}(y,t-s)|
\Phi_1(y,s) \, dy ds
&\le C(1+t)^{-1/2}\\
\end{aligned}
\end{equation}
and
\begin{equation}\label{occoneq2}
\begin{aligned}
\int_0^t\int_{-\infty}^{+\infty} |\tilde G(x,t-s;y)|
\Phi_2(y,s) \, dy ds
&\le C(\theta+\psi_1+\psi_2)(x,t),\\
\int_0^t\int_{-\infty}^{+\infty} |e_{t}(y,t-s)|
\Phi_2(y,s) \, dy ds
&\le C(1+t)^{-3/2},\\
\int_0^t\int_{-\infty}^{+\infty} |e(y,t-s)-e(y,+\infty)|
\Phi_2(y,s) \, dy ds
&\le C(1+t)^{-3/2},\\
\end{aligned}
\end{equation}
where
\begin{equation}\label{Phi}
\begin{aligned}
\Phi_1(y,s)&:= 
e^{-\eta|y|}s^{-1/2}(\theta+\psi_1+\psi_2)(y,s)
\le
Ce^{-\eta|y|/2}s^{-1/2}(1+s)^{-1},\\
\Phi_2(y,s)&:= 
e^{-\eta|y|}(1+s)^{-3/2}.\\
\end{aligned}
\end{equation}
\end{lemma}

Specifically, defining
\begin{equation}
\label{zeta2oc}
 \zeta(t):= \sup_{y, 0\le s \le t}
 \Big( |u|(\theta+\psi_1+\psi_2)^{-1}(y,t)
 + |\dot \delta (s)|(1+s)+|\delta(s)-\d|(1+s)^{1/2} \Big),
\end{equation}
noting that
$$
|Q^\d|\le \zeta^2\Psi,\quad |R^\d|\le \zeta^2(\Psi+ \Phi_1), \quad
|S^\d|\le \zeta^2\Phi_2,
$$
and applying our convolution lemmas,
we may obtain \eqref{claim} as before, yielding at once global
existence and the claimed rates of decay.
\medskip

\begin{remark}\label{coords}
It was crucial in the argument to linearize about the
limiting profile $\bar u^{\delta_{+\infty}}$ in order that
error $S$ be manageable, i.e., $\delta \to \d$ as $t\to +\infty$.
\end{remark}

\bigbreak
\section{Mixed type profiles, constant viscosity case.}\label{mixed}

We now present an alternative proof 
subsuming Lax, undercompressive, overcompressive, and even mixed
under--overcompressive cases in a single argument.
For clarity of exposition, we first restrict to the 
simpler case $B\equiv \constant$, which permits also the following 
slightly stronger result.
The general case is treated in Section \ref{gmixed}.

\begin{theorem}\label{nonlinmixed}
Let $B\equiv \constant$.  Then,
assuming (H0)--(H4), and stability condition ($\mathcal{D}$),
the profile $\bar u$ is nonlinearly phase-asymptotically
orbitally stable with respect to (not necessarly H\"older continuous)
initial perturbations $|u_0(x)|\le E_0 (1+|x|)^{-3/2}$, $E_0$ 
sufficiently small.
More precisely, \eqref{pointwise} is satisfied
for some $\delta(\cdot)$, $\delta(+\infty)$,
	where $\tilde u$ denotes the solution of \eqref{parabolic}
with initial data $\tilde u_0=\bar u+u_0$.
\end{theorem}

{\bf Proof.}
Defining $u$ as in \eqref{ocpert}, we obtain
\begin{equation}
\label{mixedperteq}
u_t-L^\d u=Q^\d(u)_x+\dot \delta (t)
(\partial \bar u^\delta/\partial \delta)_{|\d}
+ R^\d(\delta, u)_x +
S^\d(\delta, \delta_t),
\end{equation}
where 
\begin{equation}\label{mixedQ}
Q^\d =\mathcal{O}( |u|^2),
\end{equation}
\begin{equation}\label{mixedR}
R^\d=
\Big(A(\bar u^\d(x))-A(\bar u^{\delta(t)}(x))\Big)u
=\mathcal{O}( e^{-\eta|y|}|\delta-\d||u| )
\end{equation}
and $S^\d$ is as in \eqref{S}.
Defining 
\begin{equation}
\label{mixeddelta2}
 \begin{aligned}
  \delta (t)&=
\d -\int^\infty_{-\infty}e(y,t) u_0^\d(y)\,dy 
  -\int^t_0\int^{+\infty}_{-\infty} e(y,t-s)
S^\d(\delta, \delta_t)(y,s) dy ds 
\\
  &\qquad
  +\int^t_0\int^{+\infty}_{-\infty} e_{y}(y,t-s)
(Q^\d(u)+ R^\d(\delta, u))
(y,s) dy ds, \\
u_0^\d(y)&:=(\tilde u-\bar u^\d)(y),\\
  \end{aligned}
\end{equation}
we obtain 
\begin{equation}
\label{mixedu}
\begin{aligned}
  u(x,t)&= \mathcal{T}_u(u,\delta, \dot \delta,\d)(t)
:=\int^\infty_{-\infty} \tilde G(x,t;y)
\bar u_0^\d(y)\,dy \\
&\qquad -
  \int^t_0\int^\infty_{-\infty}\tilde G(x,t-s;y)
  S^\d(\delta, \delta_t))(y,s) dy \, ds \\
 &\qquad
  -\int^t_0\int^\infty_{-\infty}\tilde G_y(x,t-s;y)(Q^\d(u)+
  R^\d(\delta, u))(y,s) dy \, ds
\end{aligned}
\end{equation}
and
\begin{equation}
\label{mixeddeltadot}
 \begin{aligned}
  \dot \delta (t) &= \mathcal{T}_{\dot\delta}(u,\delta, \dot \delta,\d)(t)
:=-\int^\infty_{-\infty}e_t(y,t) 
\bar u_0^\d(y)\,dy \\
&\qquad
  - \int^t_0\int^{+\infty}_{-\infty} e_{t}(y,t-s)S^\d(\delta,\delta_t)
  (y,s)\,dy\,ds\\
 &\qquad
  +\int^t_0\int^{+\infty}_{-\infty} e_{yt}(y,t-s)(Q^\d(u)+
  R^\d(\delta, u))(y,s)\,dy\,ds. 
 \end{aligned}
\end{equation}

Defining now
\begin{equation}\label{mixeddelta}
 \begin{aligned}
  \mathcal{T}_{\delta}&(u,\delta, \dot \delta,\d)(t):=\\
&
-\int^\infty_{-\infty}
(e(y,t)-e(y,+\infty)) 
\bar u_0^\d(y)\,dy \\
 &\qquad
 -\int^t_0\int^{+\infty}_{-\infty} (e(y,t-s)-e(y,+\infty))
S^\d(\delta,\dot\delta)(y,s)\,dy\,ds \\
 &\qquad
 +\int^t_0\int^{+\infty}_{-\infty} (e_{y}(y,t-s)-e_y(y,+\infty))
(Q^\d(u)+ R^\d(\delta, u))(y,s)\,dy\,ds \\
 &\qquad
 +\int_t^{+\infty}\int^{+\infty}_{-\infty} e(y,+\infty)
S^\d(\delta,\dot\delta) (y,s)\,dy\,ds, \\
 &\qquad
 -\int_t^{+\infty}\int^{+\infty}_{-\infty} e_y(y,+\infty)
((Q^\d(u)+R^\d(\delta, u)) (y,s)\,dy\,ds, \\
\end{aligned}
\end{equation}
and
\begin{equation}\label{mixeddeltainfty}
 \begin{aligned}
  \mathcal{T}_{\delta_{+\infty}}&(u,\delta, \dot \delta,\d)(t):=\\
&
\d
-\int^\infty_{-\infty}
e(y,+\infty)) 
\bar u_0^\d(y)\,dy \\
 &\qquad
 -\int^{+\infty}_0\int^{+\infty}_{-\infty} e(y,+\infty)
S^\d(\delta,\dot\delta)(y,s)\,dy\,ds \\
 &\qquad
 +\int^{+\infty}_0\int^{+\infty}_{-\infty} e_y(y,+\infty)
(Q^\d(u)+ R^\d(\delta, u))(y,s)\,dy\,ds, \\
\end{aligned}
\end{equation}
we may express the solution of \eqref{mixeddelta2}--\eqref{mixeddeltadot} 
equivalently as the solution of the fixed-point equation
\begin{equation}\label{fixedpt}
(u,\delta,\dot\delta)=
(\mathcal{T}_{u} ,\mathcal{T}_{\delta} ,\mathcal{T}_{\dot \delta})
(u,\delta,\dot\delta,\d).
\end{equation}
in combination with
\begin{equation}\label{ift}
\mathcal{T}_{\delta_{+\infty}}
(u,\delta,\dot\delta,\d )=\d.
\end{equation}
Defining norm 
\begin{equation}
\begin{aligned}
|(f,g,h)|_\zeta&:=|f(\theta+\psi_1+\psi_2)^{-1}|_{L^\infty(x,t)}
+
|g(t)(1+t)^{1/2}|_{L^\infty(t)}\\
&\qquad +
|h(t)(1+t)|_{L^\infty(t)}
\end{aligned}
\label{zetanorm}
\end{equation}
and Banach space 
\begin{equation}
\mathcal{B}:=\{(f,g,h):\, |f,g,h|_\zeta<+\infty\},
\label{B}
\end{equation}
we find by the estimates of the previous sections that, for 
$$
|\tilde u_0-\bar u|(x)\le E_0(1+|x|)^{-3/2},
$$
 $E_0$ sufficiently small, 
$(\mathcal{T}_{u} ,\mathcal{T}_{\delta} ,\mathcal{T}_{\dot \delta},
\mathcal{T}_{\delta_{+\infty}})$
is a well-defined mapping from 
$$
B(0,r)\subset \mathcal{B}\times \mathbb{R}\to
\mathcal{B}\times \mathbb{R}
$$
for $r>0$ sufficiently small, with
\begin{equation}
|\mathcal{T}|_{\mathcal{B}\times \mathbb{R}}=
\mathcal{O}(E_0+ |\d|+ |(u,\delta,\dot\delta)|_{\zeta}^2).
\label{small}
\end{equation}

Moreover, essentially the same estimates yield that $\mathcal{T}$
is Frechet differentiable on $B(0,r)$, with
\begin{equation}
\frac{\partial
(\mathcal{T}_{u} ,\mathcal{T}_{\delta} ,\mathcal{T}_{\dot \delta},
\mathcal{T}_{\delta_{+\infty}})}
{\partial(u,\delta, \dot\delta)}
=
{\cal O}(|(u,\delta,\dot\delta|_\zeta),
\label{contractive}
\end{equation}
\begin{equation}
\frac{\partial
(\mathcal{T}_{u} ,\mathcal{T}_{\delta} ,\mathcal{T}_{\dot \delta})}
{\partial \d }
=
{\cal O}(1),
\label{offdiag}
\end{equation}
and
\begin{equation}
\frac{\partial \mathcal{T}_{\delta_{+\infty}}}
{\partial \d }
=
{\cal O}(|(u,\delta,\dot\delta|_\zeta),
\label{diag}
\end{equation}
where the final equality \eqref{diag} follows from 
\begin{equation}
\begin{aligned}
(\partial (\mathcal{T}_{\delta_{+\infty}}-\d)/\partial \d)
&=
\int^\infty_{-\infty} e(y,+\infty)) 
(\partial \bar u^{\delta}/\partial \delta)_{|\d}(y)\,dy 
+{\cal O}(|(u,\delta,\dot\delta|_\zeta)\\
& =I_\ell 
+{\cal O}(|(u,\delta,\dot\delta|_\zeta).\\
\end{aligned}
\label{approxid}
\end{equation}

In turn, relation
$$\int^\infty_{-\infty} e(y,+\infty)) 
(\partial \bar u^{\delta}/\partial \delta)_{|\d}(y)\,dy =I_\ell
$$
follows from the standard fact that
$L^\d (\partial \bar u^{\delta}/\partial \delta)_{|\d}=0$,
hence 
$$
\int_{-\infty}^{+\infty}G(x,t;y)
(\partial \bar u^{\delta}/\partial \delta)_{|\d}(y)\, dy\equiv
(\partial \bar u^{\delta}/\partial \delta)_{|\d}(x)\, 
$$
which,
together with the fact that 
$E = (\partial \bar u^{\delta}/\partial \delta)_{|\d}(x) e(y,t))$
represents the only nondecaying
part of $G(x,t;y)$ under stability criterion ($\mathcal{D}$),
yields
$$
(\partial \bar u^{\delta}/\partial \delta)_{|\d}(x)
\int_{-\infty}^{+\infty}e(x,+\infty ;y)
(\partial \bar u^{\delta}/\partial \delta)_{|\d}(y)\, dy=
(\partial \bar u^{\delta}/\partial \delta)_{|\d}(x)\, 
$$
in the limit as $t\to +\infty$.

Combining \eqref{contractive}--\eqref{diag}, we find
that, for $E_0$, $r$ sufficiently small,
$\mathcal{T}$ is contractive with respect to norm
\begin{equation}\label{scalednorm}
|(i,j,k,l)|_*:= |(i,j,k)|_\zeta + C|l|,
\end{equation}
for $C>0$ sufficiently large, with 
$|\mathcal{T}(0)|_*=\mathcal{O}(E_0)$.
Applying the Contraction Mapping Theorem, 
we find that \eqref{fixedpt}--\eqref{ift} have a 
unique solution in $\mathcal{B}$,
from which the stated decay estimates follow by definition of $|\cdot|_\zeta$.
\myqed

\begin{remark}\label{sattinger}
The above contraction mapping argument may be recognized 
as an alternative version of the Implicit Function Theorem
argument commonly used to establish orbital stability, as for example
in \cite{He}, \cite{Sat}.
\end{remark}

\bigbreak
\section{Mixed type profiles, general case.}\label{gmixed}

The constant-viscosity assumption of the previous section
made it possible to work in a weighted $L^\infty$ norm,
since we needed to gain only a single derivative in the associated
nonlinear iteration scheme.
To treat the general case, we work instead in a H\"older space
using Schauder-type smoothing estimates like those of Lemma \ref{shorttime}.

\medbreak

{\bf Proof of Theorem \ref{nonlin}, general case.}
In the general, variable-viscosity case, we may express
the perturbation $u$ of \eqref{ocpert} again as the solution of
equation \eqref{mixedperteq}, but with
$Q^\d$ and $R^\d$ now depending also on $u_x$.

To accomodate this fact, we impose on $u_0$ also H\"older continuity,
\begin{equation}
|u_0|_{C^{0+\alpha}}\le C, 
\end{equation}
and impose on $u$ the uniform bounds
\begin{equation}
\begin{aligned}
|u|_{C^{0+\alpha}}&\le C,\qquad
|u|_{C^{1}}\le Ct^{-\frac{1}{2}+\alpha},\qquad
|u|_{C^{2}}\le Ct^{-1+\alpha},\\
\end{aligned}
\label{Cbds}
\end{equation}
or equivalently $|u|_\alpha\le C$, for
\begin{equation}
\label{alphanorm}
|u|_\alpha:= \sup_{s\ge 0}
\Big( |u|_{C^{0+\alpha}} +
|u|_{C^{1}}s^{\frac{1}{2}-\alpha}+
|u|_{C^{2}}s^{1-\alpha}\Big),
\end{equation}
together with $|(u,\delta, \dot \delta)|_{\zeta_1},$
$|\d|\le E_0$ sufficiently small, for
\begin{equation}
\begin{aligned}
|(f,g,h)|_{\zeta_1}&:=|f(\theta+\psi_1+\psi_2)^{-1}|_{L^\infty(x,t)}\\
&\qquad +
|\partial_x f t^{\frac{1}{2}}(1+t)^{-\frac{1}{2}}
(\theta+\psi_1+\psi_2)^{-1}|_{L^\infty(x,t)}\\
&\qquad +
|g(t)(1+t)^{1/2}|_{L^\infty(t)}
+ |h(t)(1+t)|_{L^\infty(t)}
\end{aligned}
\label{gzetanorm}
\end{equation}
(note: now augmented with derivative bound).
Denote
\begin{equation}
\begin{aligned}
\mathcal{B}_1&:=\{(f,g,h):\, |f,g,h|_{\zeta_1}<+\infty\}\\
\mathcal{C}&:=\{(f,g,h):\, |f|_{\alpha}<+\infty\}.\\
\end{aligned}
\label{gB}
\end{equation}

We now define $\mathcal{T}_{\delta}$, $\mathcal{T}_{\dot\delta}$, 
and $\mathcal{T}_{\delta_{\infty}}$
by \eqref{mixeddeltadot}--\eqref{mixeddeltainfty}, exactly as before.
However, we define $\mathcal{T}_u$ now implicitly, as the solution of
\begin{equation}
\label{gmixedu}
\begin{aligned}
\mathcal{T}_u(u,\delta, \dot \delta,\d)(t)
&=\int^\infty_{-\infty} \tilde G(x,t;y)
\bar u_0^\d(y)\,dy \\
&\qquad -
  \int^t_0\int^\infty_{-\infty}\tilde G(x,t-s;y)
  S^\d(\delta, \delta_t))(y,s) dy \, ds \\
 &\qquad
  -\int^t_0\int^\infty_{-\infty}\tilde G_y(x,t-s;y)(Q^\d(u, 
\partial_x \mathcal{T}_u )\\
&\qquad +
  R^\d(\delta, u, \partial_x \mathcal{T}_u))(y,s) dy \, ds,
\end{aligned}
\end{equation}
or equivalently as the solution $v$ of
\begin{equation}
\label{implicitmixedperteq}
\begin{aligned}
v_t-\hat L^\d(u) v&=\hat Q^\d(u)_x+\dot \delta (t)
(\partial \bar u^\delta/\partial \delta)_{|\d}
+ \hat R^\d(\delta, u)_x \\
&\qquad + S^\d(\delta, \delta_t),
+ \hat T^\d(u,u_x, v_x),
\end{aligned}
\end{equation}
where 
\begin{equation}\label{hatlinearized}
\hat L^{\delta_*}(u)v:= B(\bar u^\d+u) v_{xx} -(A(\bar u^\d)  v)_{x}
\end{equation}
and
\begin{equation}\label{extrasource}
\begin{aligned}
\hat T^\d(u,u_x, v_x)&=\mathcal{O}((|u|+|u_x|)|v_x|),\\
\hat T^\d(u,u_x, v_x)_x&=\mathcal{O}((|u_x|+|u_{xx}|)|v_x|
+(|u_x|+|u_{xx}|)|v_{xx}|)
\end{aligned}
\end{equation}
for $|u|$ sufficiently small.

Observing that $\hat L(u)$ may be expanded as a nondivergence-form
operator with $C^{0+\alpha}$ coefficients, we may obtain
by standard Schauder or parametrix theory,
similarly as in Lemma \ref{shorttime}, 
both short-time existence for \eqref{implicitmixedperteq} 
and also the smoothing estimates
\begin{equation}
\begin{aligned}
|v|_{C^2}(t)&\le Ct^{-1+\alpha}(|u_0|_{C^{0+\alpha}} + |u|_{\alpha}),\\
\sup_{z}|v_x|(\theta+\psi_1+\psi_2)^{-1}(z,t)
&\le Ct^{-1/2}\\
&\quad \times \Big(|u|_{\zeta_1}+
\sup_{z}|u_0| (\theta+\psi_1+\psi_2)^{-1}(z,0)\Big)
\end{aligned}
\label{mixedsmoothing1}
\end{equation}
for $0\le t\le 1$ and
\begin{equation}
\begin{aligned}
|v|_{C^2}(t)&\le C(|u_0|_{C^{0+\alpha}} + |u|_{\alpha}+ |v|_{C^1}(t-1)),\\
\sup_{z}|v_x|(\theta+\psi_1+\psi_2)^{-1}(z,t)
&\le C\\
&\quad \times \Big(|u|_{\zeta_1}+
\sup_{z}|v| (\theta+\psi_1+\psi_2)^{-1}(z,t-1)\Big)
\end{aligned}
\label{mixedsmoothing2}
\end{equation}
for $t\ge 1$.
Substituting these bounds in \eqref{gmixedu}, and performing estimates
similarly as in Section \ref{mixed}, we find by the same type
of continuous induction scheme that $\mathcal{T}$
is well-defined, and
bounded from a sufficiently large ball about the origin to itself
in $(\mathcal{B}_1\cap \mathcal{C})\times \BbbR^1$.
Moreover, similar estimates yield that $\mathcal{T}$ restricted to
this same ball in $(\mathcal{B}_1\cap \mathcal{C})\times \BbbR^1$
intersected with a sufficiently small ball in $\mathcal{B_1}\times \BbbR$
is contractive in the rescaled $\zeta_1$ norm 
\begin{equation}\label{scalednorm2}
|(i,j,k,l)|_{**}:= |(i,j,k)|_{\zeta_1} + C|l|,
\end{equation}
with $|\mathcal{T}(0)|_{\mathcal{B_1}\times \BbbR}\le CE_0$.
Combining these facts, we obtain a unique fixed-point solution
in $(\mathcal{B}_1\cap \mathcal{C})\times \BbbR^1$, for which 
$$
|(u,\delta, \dot\delta)|_{\zeta}\le CE_0,
$$
yielding similarly as in Section \ref{mixed} a global solution of the
perturbation equations satisfying the claimed bounds.
We omit the details. 
\myqed

\bigskip
\section{Technical Lemmas}\label{estimates}

We complete our analysis, finally, by the proof 
of the deferred technical lemmas 
used in Sections \ref{uc}--\ref{gmixed}.
\medbreak
\subsection{Short-time theory}\label{stime}

{\bf Proof of Lemma \ref{shorttime}.}
We carry out the proof for equations \eqref{delta}--\eqref{u}.
The proof for equations \eqref{ocdeltaalt}, \eqref{ocu} 
goes similarly.
By standard Schauder or parametrix theory 
\cite{GM,F},
there exists 
a unique solution
$\tilde u\in C^{2+\alpha}(x)\cap C^{1+\frac{\alpha}{2}}(t)$
for $t>0$ sufficiently small
of the original (unshifted) equation \eqref{parabolic},
extending  so long as 
$|\tilde u|_{C^{0+\alpha}(x)}$ remains bounded and satisfiying uniform bounds 
\begin{equation}\label{layerB}
\begin{aligned}
|\tilde u|_{C^{0+\alpha}(x)}&\le C,\quad
|\tilde u|_{C^{1}(x)}\le C\Big(\frac{t}{1+t}\Big)^{-\frac{1}{2}+\alpha},\quad
|\tilde u|_{C^{2}(x)}\le C\Big(\frac{t}{1+t}\Big)^{-1+\alpha}\\
\end{aligned}
\end{equation}
depending only $\sup |\tilde u|_{C^{0+\alpha}(x)}$.

This in turn determines $\dot \delta$, hence $\delta$, 
through \eqref{deltadot}
by a straightforward contraction-mapping/continuation argument.
(Note: in establishing contractivity for small time of the righthand
side of \eqref{deltadot}, we must establish bounds of form
$$
C\int_0^t\int_{-\infty}^{+\infty} |e_y| 
(|\tilde u_x|^2 + |\tilde u||\tilde u_{xx}|) \, dy\,ds
\le 
\tilde C\int_0^t (|\tilde u_x|_{L^\infty}^2 + 
|\tilde u|_{L^\infty}|\tilde u_{xx}|_{L^\infty}) \,ds<1
$$
for $t>0$ sufficiently small,
which follow by $|e_y|_{L^1}\le C$, a consequence of Remark \ref{eboundsrmk}, 
and
integrability of $|\tilde u_x|_{L^\infty}^2 + 
|\tilde u|_{L^\infty}|\tilde u_{xx}|_{L^\infty}$, in turn
a consequence of estimates \eqref{layerB}.)
{}From \eqref{deltadot}, moreover, we easily obtain that
$\delta\in {C^{1+\frac{\alpha}{2}}}$, and uniformly bounded for short time,
so that perturbation 
$\tilde u(x+\delta(t), t)- \bar u(x)$, or equivalently its shift
$ \tilde u(x, t)- \bar u(x-\delta(t))$,
enjoys the same regularity properties as $\tilde u$.  
(Recall that $\bar u\in C^3$, as a solution of the
$C^2$ traveling-wave ODE.)
Moreover, $|\tilde u|_{C^{0+\alpha}}$ remains uniformly 
bounded so long as $|u|_{C^{0+\alpha}}$ does, and vice versa,
since their difference $\bar u(x-\delta(t))$ is uniformly
bounded in $C^{0+\alpha}$.

This verifies the first assertion.  To verify the second,
observe for fixed $t_0$, $\tau$ that $u(x, t)$ for 
$t_0-\tau\le t\le t_0$ by Duhamel's principle satisfies
$$
\begin{array}{l}
  \displaystyle{
  u(x,t)=\int^\infty_{-\infty}\tilde G(x,t-(t_0-\tau);y)u(y, t_0-\tau)\,dy } \\
  \displaystyle{\qquad
  +\int^t_{t_0-\tau} \int^\infty_{-\infty} \tilde G_y(x,t-s;y)
  (Q(u,u_x)+\dot \delta u ) (y,s)\,dy\,ds}
\end{array}
$$
and therefore
$$
\begin{array}{l}
  \displaystyle{
  u_x(x,t)=
\int^\infty_{-\infty}\tilde G_x(x,t-(t_0-\tau);y)u(y, t_0-\tau)\,dy } \\
  \displaystyle{\qquad
  +\int^t_{t_0-\tau} \int^\infty_{-\infty} \tilde G_x(x,t-s;y)
  (Q(u,u_x)_x+\dot \delta u_x ) (y,s)\,dy\,ds}.
\end{array}
$$

Combining \eqref{Q} and \eqref{Cbds}, we obtain
\begin{equation}
\begin{aligned}
  |u_x(x,t)|&\le
\int^\infty_{-\infty}|\tilde G_x(x,t-(t_0-\tau);y)||u(y, t_0-\tau)|\,dy  \\
  &\quad
  +C\int^t_{t_0-\tau} \int^\infty_{-\infty} |\tilde G_x(x,t-s;y)|\\
&\quad \times
\Big(|u_x| (s-(t_0-\tau))^{-\frac{1}{2}+\alpha}+
|u| (s-(t_0-\tau))^{-1+\alpha}\Big)(y,s)
\,dy\,ds,
\end{aligned}
\end{equation}
from which we readily obtain 
$$
|u_x\Psi^{-1}(x,t)|\le (t-(t_0-\tau))^{-\frac{1}{2}}
\sup_z |u(z,t_0-\tau)\Psi^{-1}(z,t_0-\tau)|
$$
for $t_0-\tau\le t\le t_0$ by a contraction argument based on the
convolution estimates of Lemmas \ref{hzlemma}--\ref{interaction2},
thus verifying \eqref{deriv}.  (By uniqueness, the fixed point
obtained by contraction mapping must in fact be $u$.)
See Lemma 5.1 \cite{SzZ} or Lemma 11.5, \cite{ZH} for similar arguments.

Likewise, we obtain by quite similar argument the H\"older bound
$$
|\partial_t^{\frac{\alpha}{2}}u|| \Psi^{-1}(x,t)|\le (t-(t_0-\tau))^{-\alpha}
\sup_z |u(z,t_0-\tau)\Psi^{-1}(z,t_0-\tau)|,
$$
where $\partial_t^{\frac{\alpha}{2}}u(x,t)$ denotes 
$
 \limsup_{\epsilon\to0}
|u(x,\cdot)|_{C^{0+\frac{\alpha}{2}} [t-\epsilon, t+\epsilon]}
$
with some abuse of notation.
Combining this fact with the uniform bound
$ |{\Psi_t }{\Psi}^{-1} |\le C $
obtainable by direct calculation,
we obtain the claimed continuity of $\sup_z |u\Psi^{-1}(z,\cdot)|$.
Indeed, we obtain H\"older continuity, $C^{0+ \frac{\alpha}{2}}$.
\myqed

\medbreak
\subsection{Integral estimates}\label{integral}
Throughout the analysis, we will make use of the following lemmas.

\begin{lemma}\label{hzlemma} Let $f(y)\ge 0$ be a bounded, 
nonincreasing function on $\mathbb{R}_+$, and also let $f \in L^1(\mathbb{R})$. 
Then for any $a > 0$ and $z > 0$, and for any $\omega > 1$,
\begin{equation*}
\begin{aligned}
\int_0^{+\infty} a^{1/2} &e^{-a(z-y)^2} f(y) dy
\le
(\frac{\sqrt{\pi}}{2} f(z/\omega)) \wedge  (a^{1/2} \|f\|_{L^1(\mathbb{R})}) \\
& +[(\frac{\sqrt{\pi}}{2} \|f\|_{L^\infty(\mathbb{R})}) 
\wedge (a^{1/2} \|f\|_{L^1(\mathbb{R})})] e^{-a\gamma z^2},
\end{aligned}
\end{equation*} 
for any $\gamma < (1-\frac{1}{\omega})^2$, and where $\wedge$ represents
{\it minimum}. 
\end{lemma}

{\bf Proof of Lemma \ref{hzlemma}.}  Lemma \ref{hzlemma} is 
proven as Lemma 6.3 in \cite{HZ}. \myqed

\medskip

The following two lemmas, useful in analyzing the $\theta(x,t)$ terms in 
the nonlinearity $\Psi$, 
can be proven in straightforward fashion by completing an appropriate 
square.  

\begin{lemma}\label{interaction1} For any $x$, $y$, $s$, $t$, $M_1$, $M_2$
$a$, and $b$, we have 
\begin{equation*}
\begin{aligned}
&\frac{(x-y-a(t-s))^2}{M_1(t-s)}+\frac{(y-bs)^2}{M_2s} 
=
\frac{(x-a(t-s)-bs)^2}{M_1(t-s)+M_2s} \\
& \quad +
\frac{M_1(t-s)+M_2s}{M_1M_2s(t-s)}\Big{(}y-\frac{(xM_2s-(aM_2+bM_1)(t-s)s)}{M_1(t-s)+M_2s}\Big{)}^2.
\end{aligned}
\end{equation*}
\end{lemma}

\smallskip

\begin{lemma}\label{interaction2} For any $x$, $y$, $s$, $t$, $M_1$, $M_2$
$a$, $b$, and $c$ we have 
\begin{equation*}
\begin{aligned}
&\frac{(x-\frac{a}{b}y-a(t-s))^2}{M_1(t-s)} + \frac{(y-cs)^2}{M_2s} 
=
\frac{(x-a(t-s)-c\frac{a}{b}s)^2}{M_1(t-s)+M_2(\frac{a}{b})^2s} \\
&\quad\quad +
\frac{M_1(t-s)+M_2(\frac{a}{b})^2s}{M_1M_2(\frac{a}{b})^2s(t-s)} \\
&\times \Big{(} \frac{a}{b}y- \frac{(a(\frac{a}{b})^2M_2 
  + (a(t-s)+c\frac{a}{b}s)M_1)s(t-s)-x(\frac{a}{b})^2M_2s}{M_1(t-s)+(\frac{a}{b})^2M_2s} \Big{)}^2.
\end{aligned}
\end{equation*}
\end{lemma}

{\bf Proof of Lemma \ref{iniconvolutions}.}
In each case of Lemma \ref{iniconvolutions}, we proceed in the 
case $x,y\le 0$.  The case $y\le 0 \le x$ is similar to the 
reflection estimate for $x,y \le 0$.  The case $y>0$ is entirely symmetric. 
(See Remark \ref{GFremarks} for a discussion of our terminology
regarding convection, reflection, and transmission kernels.  We 
will designate the case $\gamma=0$ as the Lax case and the case 
$\gamma=1$ as the undercompressive case.) 

For the first estimate in Lemma \ref{iniconvolutions},
we consider the integrals
\begin{equation*}
\int_{-\infty}^0 |\tilde{G} (x,t;y)| (1+|y|)^{-3/2} dy.
\end{equation*}

{\it Convection estimate.} For the convection kernel, 
according to Lemma \ref{hzlemma}, we have 
\begin{equation*}
\begin{aligned}
&\int_{-\infty}^0 t^{-1/2} e^{-\frac{(x-y-a_k^- t)^2}{Mt}} (1+|y|)^{-3/2} dy \\
&\quad \le
C\Big{(}t^{-1/2}\wedge(1+|x-a_k^- t|)^{-3/2} 
 + (1+t)^{-1/2} e^{-\frac{(x-a_k^- t)^2}{M't}} \Big{)},
\end{aligned}
\end{equation*}
where $M' > M$ can be taken as close to $M$ as we choose (by choosing
the $\gamma$ of Lemma \ref{hzlemma} sufficiently close to 1).
In the event that $x$ and $a_k^-$ have opposite signs, we have decay
of the form $(|x|+t)^{-3/2}$, which can be absorbed by the claimed
estimates.  In the event that $x$ and $a_k^-$ have the same sign, 
the terms $(1+t)^{-1/2} e^{-\frac{(x-a_k^- t)^2}{Mt}}$ are exactly 
the $\theta(x,t)$.  In order to see that the expression 
$t^{-1/2}\wedge(1+|x-a_k^-t|)^{-3/2}$ can be absorbed into the sum
$\theta + \psi_1 + \psi_2$, we first observe that for 
$|x - a_k^- t| \le C\sqrt{t}$, 
$$
t^{-1/2}\wedge(1+|x-a_k^-t|)^{-3/2}
\le
C_1 (1+t)^{-1/2} e^{-\frac{(x-a_k^- t)^2}{Mt}},
$$
for some constant $C_1$.
On the other hand, for $|x-a_k^- t| > C\sqrt{t}$, 
$$
(1+|x-a_k^-t|)^{-3/2}
\le
t^{-1/2}(1+|x-a_k^-t|)^{-1/2},
$$ 
which is sufficient in the case $|x| \le |a_1^-| t$.
Finally, for $x\in [\frac{a_k^-}{2}t,2a_k^- t]$, 
$t^{-1/2} \le C(|x|+t)^{-1/2}$, while for 
$x\notin [\frac{a_k^-}{2}t,2a_k^- t]$, there can
be only limited cancellation between $x$ and $a_k^- t$,
and we have   
$$
(1+|x-a_k^- t|)^{-3/2} \le C (1+|x|+t)^{-3/2}.
$$
For $|x| \ge |a_1^-| t$, we analyze the case $|x-a_1^- t|\le C\sqrt{t}$
as above, leaving the case $x \le a_1^- t - C\sqrt{t}$, for which 
$$
t^{-1/2} \wedge (1+|x-a_1^- t|)^{-3/2}
\le 
C (1 + |x - a_1^- t| + t^{1/2})^{-3/2}.
$$

{\it Reflection estimate.} In the case of the reflection kernel, Lemma \ref{hzlemma}
provides the estimate
\begin{equation*}
\begin{aligned}
&\int_{-\infty}^0 t^{-1/2} e^{-\frac{(x-\frac{a_j^-}{a_k^-}y-a_j^- t)^2}{Mt}} (1+|y|)^{-3/2} dy \\
&\quad \le
C\Big{(}t^{-1/2}\wedge(1+|x-a_j^- t|)^{-3/2} 
 + (1+t)^{-1/2} e^{-\frac{(x-a_j^- t)^2}{M't}} \Big{)},
\end{aligned}
\end{equation*}
which can be analyzed exactly as above.  For the transmission kernel (and $x<0$), 
we have exponential decay in $x$.  In the event that $|x|\ge \epsilon t$
for some $\epsilon > 0$, we have exponential decay in both $t$ and $x$, 
which can be subsumed.  In
the event that $|x| \le \epsilon t$, we have 
\begin{equation*}
\begin{aligned}
&\int_{-\infty}^0 t^{-1/2} e^{-\frac{(x-\frac{a_j^+}{a_k^-}y-a_j^+ t)^2}{Mt}} e^{-\eta|x|} (1+|y|)^{-3/2} dy \\
&\quad \le
C(1 + t)^{-3/2} e^{-\eta |x|},
\end{aligned}
\end{equation*}
which can be subsumed into the claimed estimates.  The cases for $x \ge 0$ 
are similar.

For the second estimate of Lemma \ref{iniconvolutions}, we consider integrals of the
form 
\begin{equation*}
\int_{-\infty}^0 |e_t (y,t)| (1+|y|)^{-3/2} dy,
\end{equation*} 
where according to Remark \ref{eboundsrmk}
$$
|e_t (y,t)| \le C t^{-1/2} \sum_{a_k^- > 0} e^{-\frac{(y+a_k^- t)^2}{Mt}}.
$$
According to Lemma \ref{hzlemma}, we can estimate this integral as, 
\begin{equation*}
\begin{aligned}
\int_{-\infty}^0 |e_t (y,t)| (1+|y|)^{-3/2} dy
&\le
C \int_{-\infty}^0 t^{-1/2} e^{-\frac{(y+a_k^- t)^2}{Mt}} (1+|y|)^{-3/2} dy \\
&\quad \le
C (1+t)^{-3/2}.
\end{aligned}
\end{equation*}
The third estimate of Lemma \ref{iniconvolutions} follows directly from the boundedness 
of $e(y,t)$.

For the final estimate of Lemma \ref{iniconvolutions}, we consider 
integrals of the form 
\begin{equation*}
\int_{-\infty}^0 |e(y,t) - e(y,+\infty)| (1+|y|)^{-3/2},
\end{equation*}
where according to Remark \ref{eboundsrmk}, 
$$
|e(y,t)-e(y,+\infty)|\le 
C\textrm{errfn }\left(\frac{|y|- at}{M\sqrt{t}}\right)
$$
for some $C$, $M$, $a>0$.
Computing directly, we have, then 
\begin{equation*}
\begin{aligned}
\int_{-\infty}^{0}
&|e(y,t)-e(y,+\infty)|(1+|y|)^{-3/2}\, dy \\
&\le
C \int_{-\infty}^{+\infty}
\textrm{errfn }\left(\frac{|y|- at}{M\sqrt{t}}\right) (1+|y|)^{-3/2}\, dy\\
&=
\frac{1}{2\pi} C\int_{-\infty}^0 \Big{(}\int_{-\infty}^{\frac{-y-at}{M\sqrt{t}}}e^{-z^2}dz\Big{)}
(1+|y|)^{-3/2} dy \\
&=
\frac{1}{2\pi} C \int_{-\infty}^{-\frac{a}{2}t} \Big{(}\int_{-\infty}^{\frac{-y-at}{M\sqrt{t}}}e^{-z^2}dz\Big{)}
(1+|y|)^{-3/2} dy \\
&+
\frac{1}{2\pi} C \int_{-\frac{a}{2}t}^0 \Big{(}\int_{-\infty}^{\frac{-y-at}{M\sqrt{t}}}e^{-z^2}dz\Big{)}
(1+|y|)^{-3/2} dy \\
&\le C_1(1+t)^{-1/2} + C_2 e^{-\frac{a^2}{8M^2}t}.
\end{aligned}
\end{equation*}

This completes the proof for $x,y \le 0$.  As discussed above, the cases in which 
either $x$ or $y$ is positive follow similarly.
\myqed

\medskip

{\bf Proof of Lemma \ref{convolutions}.} For Lemma \ref{convolutions}, the proof
of each estimate requires the analysis of several cases.  We proceed by 
carrying out detailed calculations in the most delicate cases and sufficing
to indicate the appropriate arguments in the others.  In particular, we 
will always consider the case $x,y \le 0$.  The case $y \le 0 \le x$ is similar
(though certainly not identical) to the reflection case for $x,y \le 0$. 
The estimates for $y\ge 0$ are entirely symmetric.  For the nonlinearity 
$\Psi$ we observe the inequality
$$
\Psi \le C (1+s)^{1/2} s^{-1/2} (\theta^2 + \psi_1^2 + \psi_2^2).
$$

{\it Nonlinearity }$\theta^2$.  We begin by estimating convolutions of the 
form
$$
\int_0^t \int_{-\infty}^0 |\tilde{G}_y (x,t-s;y)| 
 (1+s)^{1/2} s^{-1/2} \theta(y,s)^2 dy ds.
$$

{\it Lax convection}, $x\le 0$.  For $a_j^- < 0$, we consider
integrals of the form 
\begin{equation*}
\int_0^t \int_{-\infty}^0 (t-s)^{-1} e^{-\frac{(x-y-a_k^-(t-s))^2}{M(t-s)}}
  (1+s)^{-1/2} s^{-1/2} e^{-\frac{(y-a_j^-s)^2}{Ms}} dy ds,
\end{equation*}
where the constant $M$ arising in $\theta^2$ is larger than 
the exact constant, $L/2$ (making the term an upper estimate).
According to Lemma \ref{interaction1}, we have
\begin{equation*}
\begin{aligned}
&\int_0^t \int_{-\infty}^0 (t-s)^{-1} e^{-\frac{(x-y-a_k^-(t-s))^2}{M(t-s)}}
  (1+s)^{-1/2} s^{-1/2} e^{-\frac{(y-a_j^-s)^2}{Ms}} dy ds \\	
&\quad =
\int_0^t \int_{-\infty}^0 (t-s)^{-1} e^{-\frac{(x-a_k^- (t-s) -a_j^- s)^2}{Mt}}
  (1+s)^{-1/2} s^{-1/2} \\
& \quad \quad \times e^{-\frac{t}{Ms(t-s)}(y - \frac{xs - (a_k^- +a_j^-)s(t-s)}{Mt})^2} dy ds \\
&\quad \le
Ct^{-1/2} \int_0^t (t-s)^{-1/2} (1+s)^{-1/2} 
 e^{-\frac{(x-a_k^-(t-s)-a_j^-s)^2}{Mt}} ds.
\end{aligned}
\end{equation*}
The convection rate $a_k^-$ can be either positive or negative, 
so we divide the analysis into three cases: 
(1) $a_j^- < 0 < a_k^-$, (2) $a_k^- \le a_j^- < 0$, and 
(3) $a_j^- < a_k^- < 0$.
For the first, we observe that for $|x| \ge |a_j^-| t$, we can 
write 
$$
x-a_k^-(t-s)-a_j^-s=(x-a_j^-t)-(a_k^- - a_j^-) (t-s).
$$    
In this last expression, $x-a_j^-t$ and $-(a_k^- - a_j^-)(t-s)$ are 
both negative and cannot cancel, and we can compute 
\begin{equation*}
\begin{aligned}
&Ct^{-1/2} \int_0^t (t-s)^{-1/2} (1+s)^{-1/2} 
 e^{-\frac{(x-a_k^-(t-s)-a_j^-s)^2}{Mt}} ds \\
&\quad \le
Ct^{-1/2} \int_0^t (t-s)^{-1/2} (1+s)^{-1/2} 
 e^{-\frac{(x-a_j^-t)^2}{Mt}} ds \\
&\quad \le
C (1+t)^{-1/2} e^{-\frac{(x-a_j^-t)^2}{Mt}}.
\end{aligned}
\end{equation*}
(We observe that the seeming blow-up of $t^{-1/2}$ as $t \to 0$
is compensated for by the interval of integration $s \in [0,t]$.)
In the event that $|x| \le |a_j^-| t$, we further divide the 
analysis into cases $s\in [0,t/2]$ and $s \in [t/2,t]$.  
For
$s \in [0,t/2]$, we write
$$
x-a_k^-(t-s)-a_j^-s=(x-a_k^-t)+(a_k^- - a_j^-)s,
$$
for which $x-a_k^-t \le 0$ and $(a_k^- - a_j^-)s$ have
opposite signs and cancellation occurs.  In this way, we have 
the balance estimate
\begin{equation*}
\begin{aligned}
&(1+s)^{-1/2} e^{-\frac{(x-a_k^-(t-s)-a_j^-s)^2}{Mt}} \\
& \quad \le
C\Big{[} (1+|x-a_k^- t|)^{-1/2} e^{-\frac{(x-a_k^-(t-s)-a_j^-s)^2}{Mt}}
 + (1+s)^{-1/2} e^{-\frac{(x-a_k^-t)^2}{M't}} \Big{]}.
\end{aligned}
\end{equation*}
We observe that this kind of balance estimate is contained within the 
proof of Lemma \ref{hzlemma}, but that the lemma cannot quite be directly 
applied, because we first must recognize which variables balance.  In
this way, the current analysis is a refinement of the analysis of 
\cite{HZ}.  We also observe that it is
critical that for $s\in [0,t/2]$ the cancellation comes from $s$
(which on this interval does not yield $t$ decay), while for 
$s \in [t/2,t]$ it will be critical that the cancellation comes from 
$(t-s)$.  Indeed, at a purely technical level, this observation 
drives the analysis.  Computing directly, we now have 
\begin{equation*}
\begin{aligned}
&t^{-1/2} \int_0^{t/2} (t-s)^{-1/2} (1+s)^{-1/2} 
 e^{-\frac{(x-a_k^-(t-s)-a_j^-s)^2}{Mt}} ds \\
&\le 
C t^{-1/2} \int_0^{t/2} (t-s)^{-1/2}
 \Big{[} (1+|x-a_k^- t|)^{-1/2} e^{-\frac{(x-a_k^-(t-s)-a_j^-s)^2}{Mt}} ds \\
&\quad \quad +
(1+s)^{-1/2} e^{-\frac{(x-a_k^-t)^2}{M't}} \Big{]} ds.
\end{aligned}
\end{equation*}  
Finally, observing that in this case $|x - a_k^- t| \ge c(|x|+t)$, and that 
$t \ge |x|/|a_j^-|$ we determine an estimate by 
$$
(1+t+|x|)^{-1}, 
$$
which can be subsumed into $\psi_1$.  In the case $s \in [t/2,t]$ we argue
similarly, beginning with the relation 
$$
x-a_k^-(t-s)-a_j^-s = (x-a_j^- t) - (a_k^- - a_j^-) (t-s).
$$

For the second case ($a_j^- \le a_k^- < 0$), we again begin with 
the subcase $|x| \ge |a_j^-| t$, for which the argument for 
the case $a_j^- < 0 < a_k^-$ holds. Similarly, for the subcase 
$|x|\le |a_k^-| t$ we write 
$$
x-a_k^- (t-s) -a_j^- s = (x-a_k^- t) + (a_k^- - a_j^-)s,
$$ 
for which again there is no cancellation, and we obtain an estimate
by $\theta$.  In the event that $|a_k^-| t \le |x| \le |a_j^-| t$,
we proceed in the cases $s \in [0,t/2]$ and $s \in [t/2,t]$ exactly
as in the case $a_j^- < 0 < a_k^-$.  The third case, $a_k^- < a_j^- < 0$
follows similarly, though here we begin with the cases
$|x|\ge |a_k^-| t$ and $|x| \le |a_j^-| t$.

{\it Lax reflection}, $x \le 0$.  For $a_l^- < 0$, $a_j^- < 0$, and $a_k^- > 0$,
we compute from Lemma \ref{interaction2} (followed by direct integration)
\begin{equation*}
\begin{aligned}
&\int_0^t \int_{-\infty}^0 (t-s)^{-1} e^{-\frac{(x-\frac{a_j^-}{a_k^-}y-a_j^- (t-s))^2}{M(t-s)}}
 (1+s)^{-1/2} s^{-1/2} e^{-\frac{(y-a_l^- s)^2}{Ms}} dy ds \\
&\quad \le
Ct^{-1/2} \int_0^t (t-s)^{-1/2} (1+s)^{-1/2} 
 e^{-\frac{(x-a_j^-(t-s)-a_l^-\frac{a_j^-}{a_k^-}s)^2}{M((t-s)+(\frac{a_j^-}{a_k^-})^2)s}} ds.
\end{aligned}
\end{equation*}
In the case $|x| \ge |a_j^-|t$, we write 
$$
x - a_j^- (t-s) - a_l^- \frac{a_j^-}{a_k^-}s 
 = (x-a_j^- t) + (a_j^- - a_l \frac{a_j^-}{a_k^-})s,
$$
for which there is no cancellation and we obtain 
an estimate of the form of a diffusion wave.  We observe
that for $|a_j^-/a_k^-| > 1$ the estimate arising from this integral
will in general be broader than the diffusion wave we began with
(we will get a constant greater than $L$).
We can proceed, however, by dividing the integration over $s$ into cases, 
$s \in [0,t/C]$, $C$ sufficiently large, for which we have control
over the breadth of the diffusion wave, and $s \in [t/C,t]$ for 
which we have exponential decay in $t$, which can be subsumed by 
our estimates.  (For $x\ge C_1 t$, $C_1$ sufficiently large, we clearly
have exponential decay in both $|x|$ and $t$; hence, we can take decay
in $t$ to give similar decay in $|x|$).  

For the case $|x| \le |a_j^-| t$, we first consider the case 
$s \in [0,t/C]$, for which we write
$$
x-a_j^- (t-s) - a_l\frac{a_j^-}{a_k^-}s
 = (x-a_j^- t) + (a_j^- - a_l^- \frac{a_j^-}{a_k^-})s,
$$    
and proceed as in the Lax convection case with 
$s \in [0,t/2]$.  For $s \in [t/C,t-t/C]$, we clearly
have $(t+|x|)^{-1}$ decay, which can be subsumed, 
while for $s \in [t-t/C,t]$ we write
$$
x-a_j^- (t-s) -a_l^- \frac{a_j^-}{a_k^-}s
 = (x - a_l^- \frac{a_j^-}{a_k^-}s) - a_j^- (t-s),
$$
and proceed as in the Lax convection case with
$s \in [t/2,t]$.

{\it Lax transmission}, $x\le 0$.  For $a_l^- < 0$, $a_k^- > 0$, 
and $a_j^+ > 0$, we consider integrals of the form 
\begin{equation*}
\begin{aligned}
&\int_0^t \int_{-\infty}^0 (t-s)^{-1} e^{-\frac{(x-\frac{a_j^+}{a_k^-}y-a_j^+ t)^2}{Mt}} e^{-\eta |x|}
 (1+s)^{-1/2} s^{-1/2} e^{-\frac{(y-a_l^- s)^2}{Ms}} dy ds \\
&\quad \le
Ct^{-1/2} e^{-\eta |x|} \int_0^t (t-s)^{-1/2} (1+s)^{-1/2} 
 e^{-\frac{(x-a_j^+(t-s)-a_l^- \frac{a_j^+}{a_k^-}s)^2}{M((t-s)+(\frac{a_j^-}{a_k^-})^2)s}} ds.
\end{aligned}
\end{equation*}
In this case, for $|x| \ge \epsilon t$, and fixed $\epsilon > 0$, we 
have exponential decay in both $|x|$ and $t$, which can be subsumed
into our estimates.  In the case $|x| \le \epsilon t$, for the 
case $s \in [0,t/2]$, we write 
$$
x - a_j^+ (t-s) -a_l^- \frac{a_j^+}{a_k^-} s
 = (x-a_j^+ t) + (a_j^+ - a_l^- \frac{a_j^+}{a_k^-})s,
$$
and proceed as in the Lax convection case with $s \in [0,t/2]$,
while for $s \in [t/2,t]$ we write 
$$
x - a_j^+ (t-s) -a_l^- \frac{a_j^+}{a_k^-} s
 = (x - a_l^- \frac{a_j^+}{a_k^-}t) + (a_l^- \frac{a_j^+}{a_k^-} - a_j^+)(t-s),
$$
and proceed as in the Lax convection case with 
$s \in [t/2,t]$.

{\it Undercompressive convection}, $x \le 0$. For $a_j^- < 0$, we consider integrals
of the form
$$
\int_0^t \int_{-\infty}^0 (t-s)^{-1/2} e^{-\frac{(x-y-a_k^-(t-s))^2}{M(t-s)}} e^{-\eta |y|}
 (1+s)^{-1/2} s^{-1/2} e^{-\frac{(y-a_j^- s)^2}{Ms}} dy ds,
$$
for which we observe the inequality 
$$
e^{-\eta |y|} e^{-\frac{(y-a_j^- s)^2}{Ms}} \le Ce^{-\eta_1 |y|} e^{-\eta_2 s}.
$$
We can now employ Lemma \ref{interaction1} and proceed as in the Lax convection 
case, taking advantage of the integrability of $e^{-\eta_1 |y|}$ in $y$
and the integrability of $e^{-\eta_2 s}$ in $s$.  Estimates in the 
undercompressive reflection case and undercompressive transmission case
follow similarly.  This concludes the analysis of the 
nonlinearity $\theta^2$.

{\it Nonlinearity }$\psi_1^2$.  We estimate convolutions of the 
form
$$
\int_0^t \int_{-\infty}^0 |\tilde{G}_y (x,t-s;y)| 
 (1+s)^{1/2} s^{-1/2} \psi_1^2 dy ds.
$$

{\it Lax convection}, $x < 0$.  For $a_j^- < 0$, we consider 
integrals of the form 
\begin{equation*}
\begin{aligned}
&\int_0^t \int_{a_1^- s}^0 (t-s)^{-1} e^{-\frac{(x-y-a_k^- (t-s))^2}{M(t-s)}}
 (1+s)^{1/2} s^{-1/2} \\
& \quad \times  (1+|y|+s)^{-1} (1+|y-a_j^- s|)^{-1} dy ds.
\end{aligned}
\end{equation*}
Clearly, the primary new element here is that we no longer
have the precision afforded by Lemmas \ref{interaction1} 
and \ref{interaction2}.  We proceed instead through balance
estimates similar again to those that arise in the 
proof of Lemma \ref{hzlemma}.
We consider three subcases: (1) $a_j^- < 0 < a_k^-$, 
(2) $a_j^- \le a_k^- < 0$, and (3) $a_k^- < a_j^- < 0$, 
beginning with the first.  We begin by considering the
interval $|x| \ge |a_1^-| t$, on which we write
\begin{equation}
x - y - a_k^- (t-s)
 = (x - a_1^- t) - (y - a_1^- t) - a_k^- (t-s).
\label{balance1}
\end{equation}
For $s \in [0,t]$ and $y \in [a_1^- s,0]$, we have
$y - a_1^- t \ge 0$ and consequently each term in 
(\ref{balance1}) is negative, and we have no cancellation.
Integrating in straightfoward fashion, then, we can 
conclude an estimate by $C \theta(x,t)$.  
For $|x| \le |a_1^-| t$, we begin by writing
$$
x-y-a_k^- (t-s) = (x - a_k^- (t-s) - a_j^- s) - (y - a_j^- s),
$$
from which we see that either $x-y-a_k^- (t-s)$ is near
$x - a_k^- (t-s) - a_j^- s$ or 
$$
|y - a_j^- s| \ge \epsilon |x - a_k^- (t-s) -a_j^- s|,
$$
for some $\epsilon > 0$.  More precisely, we have the estimate
\begin{equation}
\begin{aligned}
&(1 + |y-a_j^- s|)^{-1/2} e^{-\frac{(x-y-a_k^- (t-s))^2}{M(t-s)}} \\
&\quad \le
C\Big{[}(1+|x-a_k^- (t-s)-a_j^- s|)^{-1/2} e^{-\frac{(x-y-a_k^-(t-s))^2}{M(t-s)}} \\
&\quad + (1+|y-a_j^- s|)^{-1/2} e^{-\frac{(x-a_k^-(t-s)-a_j^-s)^2}{M'(t-s)}}
 e^{-\epsilon\frac{(x-y-a_k^-(t-s))^2}{M(t-s)}} \Big{]}.
\end{aligned}
\label{balance2}
\end{equation}
For the second piece of this estimate, we have
\begin{equation*}
\begin{aligned}
&\int_0^t \int_{a_1^- s}^0 (t-s)^{-1} e^{-\frac{(x-a_k^-(t-s)-a_j^-s)^2}{M'(t-s)}}
   e^{-\epsilon\frac{(x-y-a_k^-(t-s))^2}{M(t-s)}} \\ 
&\quad \times  (1+s)^{1/2} s^{-1/2} (1+|y|+s)^{-1} (1+|y-a_j^- s|)^{-1} dyds. 
\end{aligned}
\end{equation*}
On the interval $s\in[0,t/2]$, we integrate 
$(1+|y-a_j^- s|)^{-1}$ in $y$, while on the interval
$s\in [t/2,t]$ we integrate the remaining Gaussian kernel.  We
obtain an estimate, then, by 
\begin{equation*}
\begin{aligned}
&C \int_0^{t/2} t^{-1} (1+s)^{1/2} s^{-1/2} \log (1+C's) 
 e^{-\frac{(x-a_k^-(t-s)-a_j^-s)^2}{M'(t-s)}} ds \\
&+
C \int_{t/2}^t (1+t)^{-1} (t-s)^{-1/2} e^{-\frac{(x-a_k^-(t-s)-a_j^-s)^2}{M'(t-s)}} ds,
\end{aligned}
\end{equation*}
both of which can now be analyzed by the methods of the 
Lax convection case with nonlinearity $\theta^2$.

The critical new estimate in this case is 
\begin{equation*}
\begin{aligned}
&\int_0^t \int_{a_1^- s}^0 (t-s)^{-1} e^{-\frac{(x-y-a_k^-(t-s))^2}{M(t-s)}}
 (1+|x - a_k^- (t-s) -a_j^- s|)^{-1/2} \\
&\quad \times 
(1+s)^{1/2} s^{-1/2} (1+|y|+s)^{-1} (1+|y-a_j^- s|)^{-1/2} dyds. 
\end{aligned}
\end{equation*}
For $|x| \ge |a_j^-| t$ (though $|x| \le |a_1^-| t$; we recall 
that the case $|x| \ge |a_1^-| t$ has already been considered), 
we write
$$
x - a_k^- (t-s) -a_j^- s = (x-a_j^- t) - (a_k^- - a_j^-)(t-s),
$$
for which we have no cancellation and straightfoward integration
provides an estimate by $(1+t)^{-1/2} (1+|x-a_k^- t|)^{-1/2}$.  
For $|x|\le |a_j^-| t$, we divide the analysis into the
cases $s \in [0,t/2]$ and $s \in [t/2, t]$.  In
the case $s \in [0,t/2]$, we write
$$
x - a_k^- (t-s) - a_j^- s = (x-a_k^- t) - (a_k^- - a_j^-)s,
$$
for which (observing that for $x<0$ and $a_k^- < 0$, 
$|x-a_k^- t| = |x| + |a_k^-| t$),
\begin{equation}
\begin{aligned}
&(1+|x-a_k^- (t-s) -a_j^- s|)^{-1/2} (1+|y|+s)^{-1/2} \\
&\quad \le
C \Big{[} (1+|x|+t)^{-1/2} (1+|y|+s)^{-1/2} \\
& + (1+|x-a_k^- (t-s)-a_j^- s|)^{-1/2} (1+|y|+|x|+t)^{-1/2} \Big{]}.
\label{lax2}
\end{aligned}
\end{equation}
For the second estimate in (\ref{lax2}), we integrate $(1+|y-a_j^- s|)^{-1/2}$ to find  
\begin{equation*}
\begin{aligned}
&\int_0^{t/2} \int_{a_1^- s}^0 (t-s)^{-1} e^{-\frac{(x-y-a_k^- (t-s))^2}{M(t-s)}}
 (1+s)^{1/2} s^{-1/2} (1+|y|+s)^{-1/2} \\
& \times (1+|y-a_j^- s|)^{-1/2} 
 (1+|x-a_k^- (t-s)-a_j^- s|)^{-1/2} (1+|x|+t)^{-1/2} dy ds \\
&\quad \le
C t^{-1} (1+|x|+t)^{-1/2} \int_0^{t/2} (1+s)^{1/2} ds \\
&\quad \le
C t^{-1/2} (1+|x|+t)^{-1/2}.   
\end{aligned}
\end{equation*}
The first estimate in (\ref{lax2}) can be analyzed similarly.
In the case $s \in [t/2, t]$, we have 
$$
x - a_k^- (t-s) - a_j^- s = (x - a_j^- t) - (a_k^- - a_j^-) (t-s),
$$
for which we have 
\begin{equation*}
\begin{aligned}
&(t-s)^{-1/2} (1+|x-a_k^- (t-s) -a_j^- s|)^{-1/2} \\
&\quad \le
C\Big{[} |x-a_j^- t|^{-1/2} (1+|x - a_k^- (t-s) a_j^- s|)^{-1/2} \\
&\quad \quad + (t-s)^{-1/2} (1+|x-a_j^- t|)^{-1/2} \Big{]}
\label{lax3}
\end{aligned}
\end{equation*}
For the first of these, we estimate
\begin{equation*}
\begin{aligned}
&\int_{t/2}^t \int_{a_j^- s}^0 (t-s)^{-1/2}  
  |x-a_j^- t|^{-1/2} (1+|x - a_k^- (t-s) a_j^- s|)^{-1/2} \\
&  \times e^{-\frac{(x-y-a_k^- (t-s))^2}{M(t-s)}} (1+s)^{1/2} s^{-1/2}
  (1+|y|+s)^{-1} dyds \\
& \quad \le 
C (1+t)^{-1} |x - a_j^- t|^{-1/2} \\
&\quad \times \int_{t/2}^t  
  (1+|x - a_k^- (t-s) a_j^- s|)^{-1/2} (1+s)^{1/2} s^{-1/2} ds \\
&\quad \le
C (1+t)^{-1/2} |x-a_j^- t|^{-1/2},
\end{aligned}
\end{equation*}
and similarly for the second.
We remark that the apparent blow-up at $x = a_j^- t$ is an 
artifact of the approach and can be removed by the observation
that for $|x - a_j^- t| \le C\sqrt{t}$, we can proceed by
alternative estimates to get decay of form $\theta(x,t)$.

The remaining cases $a_j^- \le a_k^- < 0$ and $a_k^- < a_j^- < 0$
follow similarly as in the Lax convection case for $\theta^2$, with 
the arguments augmented by balance estimates along the lines of 
(\ref{balance2}).

{\it Lax reflection.} For $a_l^- < 0$, $a_k^- > 0$, and $a_j^- < 0$,
we consider the convolutions
\begin{equation*}
\begin{aligned}
&\int_0^t \int_{a_1^- s}^0 (t-s)^{-1} 
 e^{-\frac{(x-\frac{a_j^-}{a_k^-}y-a_j^- (t-s))^2}{M(t-s)}} \\
& \times (1+s)^{1/2} s^{-1/2} (1+|y|+s)^{-1} (1+|y-a_l^- s|)^{-1} dy ds.
\end{aligned}
\end{equation*}
For $|x| \ge |a_j^-| t$, we write 
$$
x - \frac{a_j^-}{a_k^-} y - a_j^- (t-s)
 = (x - a_j^- t) - (\frac{a_j^-}{a_k^-}y - a_j^- s).
$$
Here, $x - a_j^- t < 0$ and $\frac{a_j^-}{a_k^-}y - a_j^- s > 0$,
and hence we do not have cancellation, and we get 
an estimate by $C \theta(x,t)$.  For $|x| \le |a_j^-|t$,
we write
$$
x-\frac{a_j^-}{a_k^-} y - a_j^- (t-s)
 = (x-a_j^- (t-s) - a_l^- \frac{a_j^-}{a_k^-}s)
 - \frac{a_j^-}{a_k^-} (y-a_l^- s),
$$ 
from which we have the inequality
\begin{equation*}
\begin{aligned}
&(1+|y-a_l^- s|)^{-1/2} e^{-\frac{(x-\frac{a_j^-}{a_k^-}y-a_j^- (t-s))^2}{M(t-s)}} \\
&\le 
C \Big{[} (1+|x-a_j^- (t-s) - a_l^- \frac{a_j^-}{a_k^-}s|)^{-1/2}
   e^{-\frac{(x-\frac{a_j^-}{a_k^-}y-a_j^- (t-s))^2}{M(t-s)}} \\
&+ 
(1+|y-a_l^- s|)^{-1/2} e^{-\frac{(x-a_j^- (t-s) -a_l^-\frac{a_j^-}{a_k^-}s)^2}{M'(t-s)}}
e^{-\epsilon\frac{(x-\frac{a_j^-}{a_k^-}y-a_j^- (t-s))^2}{M(t-s)}}
\Big{]}.
\end{aligned}
\end{equation*}
We proceed now as in the Lax convection case for $\psi_1^2$, 
writing for $s \in [0,t/2]$
$$
x-a_j^- (t-s) - a_l^- \frac{a_j^-}{a_k^-} s
 = (x-a_j^- t) + (a_j^- - a_l^- \frac{a_j^-}{a_k^-})s
$$
and for $s \in [t/2,t]$
$$
x-a_j^- (t-s) - a_l^- \frac{a_j^-}{a_k^-} s
 = (x - a_l^-\frac{a_j^-}{a_k^-}t) +(a_l^-\frac{a_j^-}{a_k^-}- a_j^-)(t-s).
$$

{\it Lax transmission.} For $a_l^- < 0$, $a_k^- > 0$, and $a_j^+ > 0$,
we consider the convolutions
\begin{equation*}
\begin{aligned}
&\int_0^t \int_{a_1^- s}^0 (t-s)^{-1} 
 e^{-\frac{(x-\frac{a_j^+}{a_k^-}y-a_j^+ (t-s))^2}{M(t-s)}} e^{-\eta |x|} \\
& \times (1+s)^{1/2} s^{-1/2} (1+|y|+s)^{-1} (1+|y-a_l^- s|)^{-1} dy ds.
\end{aligned}
\end{equation*}
In the case $|x| \ge \epsilon t$, some fixed $\epsilon > 0$, we have 
exponential decay in both $|x|$ and $t$, which can be subsumed.  In
the case $|x| \le \epsilon t$, we proceed almost exactly as in the 
Lax reflection case for $\psi_1^2$.

{\it Undercompressive convection}.  For $a_j^- < 0$, we consider the convolutions
\begin{equation*}
\begin{aligned}
&\int_0^t \int_{a_1^- s}^0 (t-s)^{-1/2} 
 e^{-\frac{(x-y-a_k^-(t-s))^2}{M(t-s)}} e^{-\eta |y|} \\
& \times (1+s)^{1/2} s^{-1/2} (1+|y|+s)^{-1} (1+|y-a_l^- s|)^{-1} dy ds.
\end{aligned}
\end{equation*}
We observe here the inequality
$$
e^{-\eta |y|} (1+|y-a_j^- s|)^{-1} 
 \le C \Big{[}e^{-\eta_1 |y|} e^{-\eta_2 s} + e^{-\eta |y|} (1+s)^{-1} \Big{]}.
$$
Integrating $e^{-\eta |y|}$ (or $e^{-\eta_1 |y|}$), 
we now proceed similarly as in 
the analysis of the Lax convection case for $\psi_1^2$.
Similarly, the undercompressive reflection and transmission
estimates follow as in the Lax reflection and transmission
esimates.  This concludes the analysis for the nonlinearity
$\psi_1^2$.

{\it Nonlinearity }$\psi_2^2$.  We estimate convolutions of the 
form
$$
\int_0^t \int_{-\infty}^{a_1^- s} |\tilde{G}_y (x,t-s;y)| 
 (1+s)^{1/2} s^{-1/2} \psi_2(y,s)^2 dy ds.
$$

{\it Lax convection}. For $a_1^- < 0$, we consider convolutions of the form 
\begin{equation*}
\begin{aligned}
&\int_0^t \int_{-\infty}^{a_1^- s} (t-s)^{-1} 
 e^{-\frac{(x-y-a_k^-(t-s))^2}{M(t-s)}}  \\
& \quad \times (1+s)^{1/2} s^{-1/2} (1+|y-a_1^- s| + s^{1/2})^{-3} dy ds.
\end{aligned}
\end{equation*}
We write 
$$
x-y-a_k^- (t-s) 
 = (x - a_k^- (t-s) -a_1^- s) - (y-a_1^- s),
$$
for which 
\begin{equation*}
\begin{aligned}
&(1+|y-a_1^- s|+s^{1/2})^{-3/2} e^{-\frac{(x-y-a_k^- (t-s))^2}{M(t-s)}} \\
&\quad \le
C \Big{[} (1+|x-a_k^-(t-s) - a_1^- s|+s^{1/2})^{-3/2} e^{-\frac{(x-y-a_k^- (t-s))^2}{M(t-s)}} \\
&\quad +
(1+|y-a_1^- s|+s^{1/2})^{-3/2} e^{-\frac{(x-a_k^- (t-s)-a_1^- s)^2}{M' (t-s)}}
e^{-\epsilon \frac{(x-y-a_k^- (t-s))^2}{M(t-s)}} \Big{]}.
\end{aligned}
\end{equation*}
For $|x| \ge |a_1^-| t$, we write
$$
x-a_k^- (t-s) - a_j^- s
 = (x-a_1^- t) - (a_k^- - a_1^-) (t-s),
$$
for which there is no cancellation and integration gives an estimate by 
$$
C_1 (1+t)^{-1}\log(1+t) e^{-\frac{(x-a_1^- t)^2}{M't}} 
 + C_2 (1 + t)^{-1/4} (1+|x-a_1^- t| + t^{1/2})^{-3/2}.
$$
For $|x| \le |a_1^-| t$, we consider only the case $a_k^- < 0$.
The case $a_k^- > 0$ is similar.  We first consider the 
additional subcase $|x| \le |a_k^-| t$, for which we write
$$
x - a_k^- (t-s) - a_1^- s
 = (x - a_k^- t) + (a_k^- - a_1^-) s.
$$
Observing that there is no cancellation between these
terms, we can proceed exactly as in the case $|x| \ge |a_1^-| t$
to get an estimate by $C (\theta (x,t) + \psi_2)$.  
For the case $|a_k^-| t \le |x| \le |a_j^-| t$,
we divide the analysis into subcases $s \in [0,t/2]$ and
$s \in [t/2,t]$.  For $s \in [0,t/2]$, we write
$$
x - a_k^- (t-s) -a_1^- s
 = (x-a_k^- t) + (a_k^- - a_1^-)s,
$$
for which we have
\begin{equation}
\begin{aligned}
&(1 + |x-a_k^- (t-s) -a_1^- s| + s^{1/2})^{-3/2} \\
&\quad \le
C\Big{[}(1 + |x-a_k^- t| + s^{1/2})^{-3/2} \\ 
&\quad + (1 + |x-a_k^- (t-s) -a_1^- s| + |x-a_k^- t|^{1/2}  + s^{1/2})^{-3/2} \Big{]}.
\end{aligned}
\label{est2}
\end{equation}
For the first of the estimates, integrating $(1 + |y-a_1^- s| + s^{1/2})^{-3/2}$,
we estimate
\begin{equation*}
\begin{aligned}
&\int_0^{t/2} \int_{-\infty}^{a_1^- s} (t-s)^{-1} (1 + |x - a_k^- (t-s) - a_1^- s| + s^{1/2})^{-3/2}
 e^{-\frac{(x-y-a_k^- (t-s))^2}{M(t-s)}} \\
& \quad \times
(1+s)^{1/2} s^{-1/2} (1 + |y - a_1^- s| + s^{1/2})^{-3/2} dyds \\
&\quad \le
C t^{-1} \int_0^{t/2} (1+|x-a_k^- t| + s^{1/2})^{-3/2} (1+s)^{1/2} s^{-1/2} (1+s^{1/2})^{-1/2} ds \\
&\quad \le
C (1 + t)^{-3/4} (1+|x-a_k^- t|)^{-1/2}.
\end{aligned}
\end{equation*}
For $s \in [t/2,t]$, we write
$$
x - a_k^- (t-s) -a_1^- s
 = (x-a_1^- t) - (a_k^- - a_1^-)(t-s),
$$
for which we have 
\begin{equation*}
\begin{aligned}
&(t-s)^{-1/2} \Big{(}1 + |x - a_k^- (t-s) - a_1^- s| + s^{1/2} \Big{)}^{-3/2} \\
&\quad \le
C \Big{[} |x - a_1^- t|^{-1/2} \Big{(} 1 + |x - a_k^- (t-s) -a_1^- s| + s^{1/2} \Big{)}^{-3/2} \\
& \quad +
(t-s)^{-1/2} \Big{(}1 + |x-a_1^- t| + s^{1/2} \Big{)}^{-3/2} \Big{]}.
\end{aligned}
\end{equation*}
For the first, integrating over the Gaussian kernel, we estimate
\begin{equation*}
\begin{aligned}
&\int_{t/2}^t \int_{-\infty}^{a_1^- s} (t-s)^{-1} (1 + |x - a_k^- (t-s) - a_1^- s| + s^{1/2})^{-3/2}
 e^{-\frac{(x-y-a_k^- (t-s))^2}{M(t-s)}} \\
& \quad \times
(1+s)^{1/2} s^{-1/2} (1 + |y - a_1^- s| + s^{1/2})^{-3/2} dyds \\
&\quad \le
C (1+t^{1/2})^{-3/2} |x - a_1^- t|^{-1/2} \\
&\quad \times 
 \int_{t/2}^t (1+|x-a_k^- (t-s) - a_1^- s| + t^{1/2})^{-3/2} (1+s)^{1/2} s^{-1/2} ds \\
&\quad \le
C (1 + t)^{-1/2} (1+|x-a_k^- t|)^{-1/2}.
\end{aligned}
\end{equation*}
The second estimate in (\ref{est2}) can be analyzed similarly.

The Lax reflection and Lax transmission estimates follow similarly.  Finally, 
for the 
undercompressive estimates in the case of nonlinearity $\psi_2^2$, we observe
that for $y \in (-\infty,a_1^- s]$ exponential $|y|$ decay yields exponential
$s$ decay, and the estimates follow in straightforward fashion.  This concludes
the analysis for the nonlinearity $\psi_2^2$ and consequently for the 
nonlinearity 
$$
(1+s)^{1/2} s^{-1/2} (\theta + \psi_1 + \psi_2)^2.
$$

{\it Nonlinearity} $(1+s)^{-1} (\theta + \psi_1 + \psi_2)$.  We consider
convolutions of the form 
$$
\int_0^t \int_{-\infty}^0 |\tilde{G}_y (x,t-s;y)| (1+s)^{-1} (\theta + \psi_1 + \psi_2) dyds. 
$$

{\it Lax convection}.  For $a_j^- < 0$, we consider integrals of the form
$$
\int_0^t \int_{-\infty}^0 (t-s)^{-1} e^{-\frac{(x-y-a_k^- (t-s))^2}{M(t-s)}}
 (1+s)^{-1} s^{-1/2} e^{-\frac{(y-a_j^- s)^2}{Ls}} dy ds.
$$
We observe that this integral is better than the one analyzed in the Lax
convection case with nonlinearity $\theta^2$, except that the constant 
$L$ in the diffusion kernel must now be kept.  (In general, our balance 
estimates increase the size of this constant, see especially Lemma \ref{hzlemma}.)
According to Lemma \ref{interaction1}, we can estimate this integral by 
$$
C t^{-1/2} \int_0^t (t-s)^{-1/2} (1+s)^{-1} 
 e^{-\frac{1}{M(t-s)+Ls}(x-a_k^- (t-s) - a_j^- s)^2} ds. 
$$   
For $|x| \ge |a_j^-| t$ and the case $a_j^- \le a_k^-$, we 
write 
$$
x - a_k^- (t-s) -a_j^- s
 = (x - a_j^- t) - (a_k^- - a_j^-) (t-s),
$$
for which there is no cancellation and integration yields 
an estimate by 
$$
(1+t)^{-1/2} e^{-\frac{(x-a_j^- t)^2}{M't}},
$$
where $M < M' < L$.
For $|x| \le |a_j^-| t$, we divide the analysis into the 
subcases $s \in [0,t/C_1]$ and $s \in [t/C_1,t]$ for some constant 
$C_1$ sufficiently large. In the event that 
$s \in [t/C_1 , t]$, we integrate $(t-s)^{-1/2}$
to obtain the estimate
$$
C (1+t)^{-1}, 
$$ 
which can be subsumed.  On the other hand if $s \in [0,t/C_1]$,
we observe that by choice of $C_1$ the divisor
$M(t-s) + Ls$ can be kept as close to $Mt$ as we like.  
Since $M < L$ we can use this observation to recover the 
expected estimate.  Otherwise, the analysis proceeds exactly 
as in the Lax convection case for nonlinearity $\theta^2$.

The Lax reflection and Lax transmission estimates for nonlinearity 
$(1+t)^{-1} \theta$ follow from the Lax
convection argument and Lemma \ref{interaction2}.  The
undercompressive estimates follow similarly.

In the case of the nonlinearity $(1+s)^{-1} \psi_1$, we observe 
that for $y \in [a_1^- s, 0]$, 
$$
(1 + s)^{-1} \le C (1 + |y| + s)^{-1}, 
$$
and so 
$$
(1+s)^{-1} \psi_1 \le C \psi_1^2.
$$
Hence the convolution estimates for his nonlinearity follow from 
those for the nonlinearity $\psi_1^2$.  Finally, 
estimates for the nonlinearity $(1+s)^{-1} \psi_2$  are straightfoward.
This concludes the analysis of the first estimates of Lemma \ref{convolutions}.

{\it Excited term estimates}.  The remaining estimates of Lemma
\ref{convolutions} regard the {\it excited} terms $e(y,t)$.
For the integral 
$$
\int_0^t \int_{-\infty}^0 |e_{yt} (y,t-s)| \Psi(y,s) dyds,
$$
we have according to Remark \ref{eboundsrmk}
$$
|e_{yt}(y,t)| \le C(t^{-1} + \gamma t^{-1/2} e^{-\eta |y|})
 \sum_{a_k^- > 0} e^{-\frac{(y+a_k^- t)^2}{Mt}}.
$$
We observe that these estimates correspond precisely with the Lax and 
undercompressive convection kernels with $x=0$.  In this
way, we immediately obtain an estimate by
$$
C (\theta + \psi_1 + \psi_2) (0,t) \le C_1 (1+t)^{-1},
$$ 
for some constant $C_1$.
For the integral 
$$
\int_0^{\infty} \int_{-\infty}^0 |e_y (y,+\infty)| \Psi(y,s) dy ds,
$$
we have from Remark \ref{eboundsrmk} 
$$
|e_y (y,+\infty)| \le C \gamma e^{-\eta |y|},
$$
where as usual $\gamma$ is 1 for undercompressive profiles and 
0 otherwise.  For the nonlinearity $(1+s)^{1/2} s^{-1/2} \theta^2$, 
we observe the inequality
$$
e^{-\eta |y|} e^{-\frac{(y-a_j^- s)^2}{Ms}} 
\le C e^{-\eta_1 |y|} e^{-\eta_2 s},
$$
for some fixed $\eta_1 > 0$ and $\eta_2 > 0$.  
We estimate
\begin{equation*}
\int_0^{\infty} \int_{-\infty}^0 e^{-\eta_1 |y|} e^{-\eta_2 s} (1+s)^{-1/2} s^{-1/2} dyds \le C,
\end{equation*}
by the integrability of $e^{-\eta_1 |y|}$ and $e^{-\eta_2 s}$. 
Similarly, for the  nonlinearity $(1+s)^{1/2} s^{-1/2} \psi_1^2$, we use 
the inequality 
$$
e^{-\eta |y|} (1+|y-a_j^- s|)^{-1} 
 \le C \Big{[}e^{-\eta_1 |y|} e^{-\eta_2 s} + e^{-\eta |y|} (1+s)^{-1} \Big{]},
$$ 
from which the required estimate follows from the integrability of $e^{-\eta_1 |y|}$
and of $(1+s)^{-2}$.  For the nonlinearity $(1+s)^{1/2} s^{-1/2} \psi_2^2$,
we have $y \in (-\infty, a_1^- s]$, for which exponential $|y|$ decay yields 
exponential $s$ decay and the estimate is immediate.
For the nonlinearity $(1+s)^{-1} (\theta+\psi_1+\psi_2)$, we observe the 
inequality
\begin{equation*}
e^{-\eta |y|} (\theta+\psi_1+\psi_2) (y,s)
\le C e^{-\eta_1 |y|} (1+s)^{-1}.
\end{equation*}
We have, then an estimate by 
\begin{equation*}
\begin{aligned}
C& \int_0^{\infty} \int_{-\infty}^0 e^{-\eta_1 |y|} (1+s)^{-2} dy ds \\ 
&\quad \le
C_1 \int_0^\infty (1+s)^{-2} \le C_2.
\end{aligned}
\end{equation*}

For the integral 
$$
\int_0^t \int_{-\infty}^0 |e_y (y,t-s) - e_y(y,+\infty)| \Psi (y,s) dy ds,
$$
we have from Remark \ref{eboundsrmk}
$$
|e_y (y,t) - e_y(y,+\infty)| \le C t^{-1/2} \sum_{a_k^- > 0} e^{-\frac{(y+a_k^- t)^2}{Mt}}.
$$
Observing that our estimate on $|e_y (y,t) - e_y(y,+\infty)|$ takes the form 
of the Lax convection kernel multiplied by $(t-s)^{1/2}$, we proceed exactly as 
there to determine an estimate by $C (1+t)^{-1/2}$.

For the integral 
$$
\int_t^{\infty} |e_y (y,t-s)| \Psi (y,s) dy ds,
$$ 
we have from Remark \ref{eboundsrmk}
\begin{equation*}
\begin{aligned}
&|e_y (y,t)| \le C t^{-1/2} \sum_{a_k^- > 0} e^{-\frac{(y+a_k^- t)^2}{Mt}} \\
&\quad +
C \gamma e^{-\eta |y|} 
 \Big{(}\textrm{errfn} (\frac{y+a_k^- t}{\sqrt{4 \beta_k^- t}}) 
   - \textrm{errfn} (\frac{y-a_k^- t}{\sqrt{4 \beta_k^- t}}) \Big{)}.
\end{aligned}
\end{equation*}
For the first, we can proceed simiarly as in the Lax convection 
estimates for $x=0$.  For the second, we observe the inequality
$$
e^{-\eta |y|} \Psi (y,s) \le C e^{-\eta_1 |y|} (1+s)^{-3/2} s^{-1/2}.
$$
We estimate, then
\begin{equation*}
\begin{aligned}
&\int_t^{\infty} \int_{-\infty}^0 e^{-\eta_1 |y|} (1+s)^{-3/2} s^{-1/2} dy ds \\
&\quad \le
C (1+t)^{-1/2} \int_t^{\infty} (1+s)^{-1} s^{-1/2} ds \le C' (1+t)^{-1/2}.
\end{aligned}
\end{equation*}

This completes the proof of Lemma \ref{convolutions}

\myqed

{\bf Proof of Lemma \ref{occonvolutions}.}
Observing that the first three estimates of Lemma \ref{occonvolutions}
are significantly less precise than the second three, and moreover 
can be established by similar methods, we begin with the fourth, 
$$
\int_0^t \int_{-\infty}^0 |\tilde{G}_y (x,t-s;y)| \Phi_2 (y,s) dy ds,
$$
where 
$$
\Phi_2 (y,s) \le C e^{-\eta |y|} (1+s)^{-3/2}.
$$

{\it Convection estimate}. We consider integrals of the form 
$$
\int_0^t \int_{-\infty}^0 (t-s)^{-1/2} e^{-\frac{(x-y-a_k^- (t-s))^2}{M(t-s)}}
 e^{-\eta |y|} (1+s)^{-3/2} dy ds.
$$
In order to make use of the localization due to the term $e^{-\eta |y|}$, 
we use the inequality
\begin{equation}
\begin{aligned}
&e^{-\frac{(x-y-a_k^- (t-s))^2}{M(t-s)}} e^{-\eta |y|}
 \le C \Big{[} e^{-\frac{(x-a_k^- (t-s))^2}{M'(t-s)}} e^{-\eta |y|} \\
&\quad +
e^{-\frac{(x-y-a_k^-(t-s))^2}{M(t-s)}} e^{-\eta_1 |y|} e^{-\eta_2 |x-a_k^- (t-s)|} \Big{]}. 
\end{aligned}
\label{est1}
\end{equation}
For the first of these two estimates, we have the integral
\begin{equation*}
\int_0^t \int_{-\infty}^0 (t-s)^{-1/2} e^{-\frac{(x-a_k^- (t-s))^2}{M'(t-s)}}
 e^{-\eta |y|} (1+s)^{-3/2} dy ds.
\end{equation*}
In the case $|x| \ge |a_1^-| t$, we have
\begin{equation*}
\begin{aligned}
x-a_k^- (t-s) &=
 (x - a_1^- t) - (a_k^- - a_1^-) t + a_k^- s \\
& \le (x - a_1^- t) \le 0.
\end{aligned}
\end{equation*}
We have, then, an estimate by 
\begin{equation*}
\begin{aligned}
&C_1 t^{-1/2} e^{-\frac{(x-a_1^- t)^2}{M't}} \int_0^{t/2} (1+s)^{-3/2} ds \\
&\quad + C_2 (1+t)^{-3/2} e^{-\frac{(x-a_1^- t)^2}{M't}} \int_{t/2}^t (t-s)^{-1/2} ds \\
&\quad \le
C (1+t)^{-1/2} e^{-\frac{(x-a_1^- t)^2}{M't}}.
\end{aligned}
\end{equation*}
For $|x| \le |a_1^-| t$, we write 
\begin{equation*}
x - a_k^- (t-s) = (x - a_k^- t) + a_k^- s.
\end{equation*}
In the event that $a_k^- > 0$, we have 
$|x - a_k^- t| \ge c (|x| + t)$, so that 
\begin{equation*}
(1+s)^{-3/2} e^{-\frac{(x-a_k^- (t-s))^2}{M' (t-s)}} 
 \le C (|x| + t)^{-3/2}.
\end{equation*}
In this case, we have the estimate 
\begin{equation*}
C (|x| + t)^{-3/2} \int_0^t (t-s)^{-1/2} ds
 \le C (|x| + t)^{-1}. 
\end{equation*}
For $a_k^- < 0$, we first consider the case $|x| \ge |a_k^-| t$,
for which we again write
\begin{equation*}
x - a_k^- (t-s) = (x - a_k^- t) + a_k^- s.
\end{equation*}
Observing that there is no cancellation between these terms, we can 
conclude an estimate by $C \theta(x,t)$ as in the case $|x|\ge |a_1^-| t$
above.  In the event $|x| \le |a_k^-| t$, we integrate (for $\epsilon > 0$
sufficiently small)
\begin{equation*}
\begin{aligned}
&\int_0^t \int_{-\infty}^0 (t-s)^{-1/2} e^{-\frac{(x-a_k^- (t-s))^2}{M'(t-s)}}
 e^{-\eta |y|} (1+s)^{-3/2} dy ds \\
&\quad \le
C_1 t^{-1/2} \int_0^{\epsilon |x - a_k^- t|} e^{-\frac{(x-a_k^- t)^2}{M't}} (1+s)^{-3/2} ds \\
& \quad +
C_2 t^{-1/2} \int_{\epsilon |x - a_k^- t|}^{C|x-a_k^- t|\wedge t/2} (1+s)^{-3/2} ds \\
& \quad +
C_3 t^{-1/2} \int_{C|x-a_k^- t|\wedge t/2}^{t/2} e^{-\frac{(x-a_k^- t)^2}{M't}} (1+s)^{-3/2} ds \\
& \quad +
C_4 (1+t)^{-3/2} \int_{t/2}^t (t-s)^{-1/2} \\
& \quad \le
C_1 (1+t)^{-1/2} e^{-\frac{(x-a_k^- t)^2}{M' t}} + C_2 (1+t)^{-1/2} (1+|x-a_k^- t|)^{-1/2}.
\end{aligned}
\end{equation*}
For integration over the second estimate in (\ref{est1}), we proceed similarly.

In this case the analyses of the reflection estimate and 
the transmission estimate do not differ significantly 
from the analysis of the convection estimate.

{\it Overcompressive excited estimates.} For the integral
\begin{equation*}
\int_0^t \int_{-\infty}^0 |e_t (y,t-s)| e^{-\eta |y|} (1+s)^{-3/2} dy ds,
\end{equation*}
we have an estimate by 
\begin{equation*}
C \int_0^t \int_{-\infty}^0 (t-s)^{-1/2} e^{-\frac{(y+a_k^- (t-s))^2}{M(t-s)}}
e^{-\eta |y|} (1+s)^{-3/2} dy ds,
\end{equation*}
for which we observe the inequality
\begin{equation*}
e^{-\frac{(y+a_k^- (t-s))^2}{M(t-s)}} e^{-\eta |y|}
\le C \Big{[} e^{-\eta_1 |y|} e^{-\eta_2 (t-s)} \Big{]}. 
\end{equation*}
We have, then 
\begin{equation*}
\begin{aligned}
&\int_0^t \int_{-\infty}^0 (t-s)^{-1/2} e^{-\eta_1 |y|} e^{-\eta_2 (t-s)} (1+s)^{-3/2} dy ds \\
&\quad \le C_1 e^{-\frac{\eta_2}{2}t} \int_0^{t/2} (1+s)^{-3/2} ds \\
&\quad + C_2 (1+t)^{-3/2} \int_{t/2}^t (t-s)^{-1/2} e^{-\eta_2 (t-s)} ds \\
&\quad \le
C (1+t)^{-3/2}.
\end{aligned}
\end{equation*}

For the final integral
\begin{equation*}
\int_0^t \int_{-\infty}^0 |e(y,t-s) - e(y,+\infty)| e^{-\eta |y|} (1+s)^{-3/2} dy ds,
\end{equation*}
we have an estimate by 
\begin{equation*}
C \int_0^t \int_{-\infty}^0 \textrm{errfn} (\frac{|y|-a(t-s)}{M\sqrt{t-s}}) e^{-\eta |y|} (1+s)^{-3/2} dy ds,
\end{equation*}
for some constants $M > 0$, $a > 0$.  In this case, we observe the 
inequality
\begin{equation*}
\textrm{errfn} (\frac{|y|-a(t-s)}{M\sqrt{t-s}}) e^{-\eta |y|}
\le C e^{-\eta_1 |y|} e^{-\eta_2 (t-s)},
\end{equation*}
from which the estimate follows as above.

The remaining cases of Lemma \ref{occonvolutions} can be
analyzed similarly. \myqed
 
\bigbreak

{\it Acknowledgements.}
Thanks to Benjamin Texier for helpful discussions 
on the details of the fixed-point iteration scheme
and to Dan Marchesin for pointing out the interest
of multi-phase profiles with one endstate lying on
the degenerate viscosity saturation boundary (Remark \ref{strict}).
Research of P.H. was supported in part by the National Science Foundation 
under Grant No. DMS--0230003.
Research of K.Z. was supported in part by the National Science Foundation 
under Grants No. DMS-0070765 and DMS-0300487.

\noindent Peter HOWARD \\
\noindent Department of Mathematics \\
\noindent Texas A\&M University \\
\noindent College Station, TX 77843 \\
\noindent phoward@math.tamu.edu \\
\noindent and \\
\noindent Kevin ZUMBRUN \\
\noindent Department of Mathematics \\
\noindent Indiana University \\
\noindent Bloomington, IN  47405-4301\\
\noindent kzumbrun@indiana.edu \\

\end{document}